\documentclass[12pt]{article}

\usepackage{amsmath}
\usepackage{amssymb}
\usepackage{amsfonts}
\usepackage{amsthm}
\usepackage{epsf}
\usepackage{graphics}
\usepackage{graphpap}

\def\Cmpx{{\mathbb{C}}}
\def\Real{{\mathbb{R}}}
\def\Intg{{\mathbb{Z}}}
\def\natnum{{\mathbb{N}}}
\def\ratnum{{\mathbb{Q}}}
\def\cnj{\overline}
\def\union{{\cup}}

\def\tfrac#1#2{{\textstyle{\frac{#1}{#2}}}}

\def\scrhalf{{\raisebox{.3ex}%
{$\scriptstyle 1$}\!/\!\raisebox{-.3ex}{$\scriptstyle 2$}}}

\def\wt{\widetilde}
\def\wh{\widehat}

\def\spanrm{{\mathrm{span}}}
\def\innerprod(#1){{\langle #1 \rangle}}
\def\norm#1{\|#1\|}
\def\rvec(#1,#2){{\big|{}^{#1}_{#2} \big\rangle}} 
\def\lvec(#1,#2){{\big\langle{}^{#1}_{#2}\big|}} 
\def\Rvec(#1){{| #1 \rangle}} 
\def\Lvec(#1){{#1|}} 
\def\PB(#1,#2){{\{#1,#2\}}}

\newtheorem{theorem}{Theorem}
\newtheorem{lemma}[theorem]{Lemma}

\newenvironment{jglist}
{\begin{list}{$\bullet$}{\setlength{\rightmargin}{\leftmargin}}}
{\end{list}}

\def\Calg{{\cal A}^C}
\def\Nalg{{\cal A}^N}
\def\Ralg{{\cal A}^R}
\def\NCalg{{\cal A}^{NC}}
\def\Alg{{\cal A}}
\def\Balg{{\cal B}}
\def\Palg{{\Psi}}

\def\Xdom{{\cal C}_{\Real}}
\def\NXdom{{\cal C}_{\Real^2}}

\def\Pmod#1{{\Psi}_{#1}}
\def\man{{\cal M}}
\def\ad#1{{\mathrm{ad}_{#1}}}
\def\Gspace{{\cal G}}
\def\S{{\cal S}}
\def\T{{\cal T}}
\def\piT{\pi_{{}_\T}}
\def\SS(#1,#2){{\sigma^{#1}_{#2}}} 
\def\smax#1{{}_{#1}^{\uparrow}}
\def\smin#1{{}_{#1}^{\downarrow}}
\def\topI#1{{\mathbf I}_{#1}}
\def\Aop{A}
\def\Bop{B}
\def\Hamil{\cal H}

\newcounter{nn}

\title{Topology Change and Vector Modules on Noncommutative Surfaces
of Rotation}
\author{{\bf Jonathan Gratus}\thanks{email: jg@luna.ph.lancs.ac.uk}\\
\small Physics Department, Lancaster University, Lancaster LA1 4YB}

\begin{document}

\maketitle

\begin{abstract}

A non associative, noncommutative algebra is defined that may be
interpreted as a set of vector modules over a noncommutative surface
of rotation. Two of these vector modules are identified with the
analogues of the tangent and cotangent spaces in noncommutative
geometry, via the definition of an exterior derivative. This
derivative reduces to the standard exterior derivative in the
commutative limit.  The representation of this algebra is used to
investigate a simple topology change where two connected compact
surfaces coalesce to form one such surface.
\end{abstract}

\vspace{2em}

\noindent{\bf PACS:} 

\begin{tabular}{lcl}
Primary & 03.65.Fd & Algebraic Methods (Quantum mechanics) 
\\
Secondary & 02.40.Pc & General topology 
\\
& 04.60.Nc & Lattice and discrete methods (Quantum Gravity)
\end{tabular}

%\tableofcontents

\newpage

%%%%%%%%%%%%%%%%%%%%%%%%%%%%%%%%%%%%%%%%%%%%%%%%%%%%%%%%%%%%%%%%%%%%%%

\section{Introduction}
\label{ch_intr}

Noncommutative geometry has been proposed by many people
\cite{Madore_book}\cite{Connes} as a candidate for the formulation of
quantum gravity, principally because it combines the noncommutative
structure of quantum mechanics with the geometrical structure of
general relativity. Noncommutative geometry has also been proposed
\cite{DeWit1} in the membrane picture of particle mechanics.

There are two main tasks in noncommutative geometry: The first is to
find a one-parameter set of algebras $\Alg(\varepsilon)$ which are
noncommutative for $\varepsilon\ne0$ and commutative when
$\varepsilon=0$. We require that we can embed $\Alg(\varepsilon=0)$
into $C(\man,\Cmpx)$ the commutative algebra of complex valued
function on a manifold $\man$.  Here $\varepsilon$ plays the r\^ole of
$\hbar$ in quantum mechanics. This is similar to the quantisation of
Poisson manifolds and investigations of $\star$-products.  We call
$\Alg(\varepsilon)$ the noncommutative version of $\man$, although it
is the algebra not the manifold which fails to commute.  Thus if
$\man$ is the sphere $\Alg(\varepsilon)$ is called the noncommutative
sphere \cite{Gratus5}. In this article $\man$ will be a surface of
rotation. For a function $\rho:\Real\mapsto\Real$ we generate a
surface of rotation by rotating the section of the graph $y=\rho(x)$
which is below the $x$-axis, about the $x$-axis. Thus the
corresponding algebra $\Alg(\varepsilon)$ is called the noncommutative
surfaces of rotation (from now on this phrase is abbreviated to
{\bf NCSR}). These are defined in section \ref{ch_rev} and were first
described in \cite{Gratus7}, together with some simple facts such as
their Poisson structure and representations.

The second task is to write down the objects studied in differential
geometry, such as tangent bundles, cotangent bundles, exterior
algebras, metric tensors, connections and curvature, in terms of
elements of the algebra $C(\man,\Cmpx)$ and then extend these
definitions for $\Alg(\varepsilon)$.

There are two key properties required of tangent vector
fields. Firstly that they should be derivatives, i.e. follow Leibniz
rule, and secondly that they should form a module over the algebra of
functions. (That is one can multiply a vector with a scalar to give
another vector.)  It turns out that for noncommutative geometry these
two properties are incompatible, and one must choose either to have
vectors which are derivatives, or vectors which form a module.

The standard method is to choose vectors which form derivatives
\cite{Madore_book}. If the underlying algebra $\Alg$ is a matrix
algebra then it is easy to show that all such vectors are inner. That
is if $\xi:\Alg\mapsto\Alg$ such that $\xi(fg)=\xi(f)g+f\xi(g)$ for
all $f,g\in\Alg$ then there exists $h\in\Alg$ such that $\xi=\ad{h}$
where $\ad{h}f=[h,f]=hf-fh$.  Clearly if $\xi$ is inner then $f\xi$ is
not inner. In section \ref{ch_inn} we show that the same is true for
the algebra of functions on a NCSR.

In \cite{Gratus6} the author gives an alternative method of defining
tangent vectors on the noncommutative sphere. These vectors do form a
(one sided) module over the noncommutative sphere but are derivatives
only in the commutative limit. That is $\xi(fg)=\xi(f)g+f\xi(g)+
O(\varepsilon)$.  This article may be seen as the result of giving a
NCSR a vector bundle structure similar to that defined for the sphere
in \cite{Gratus6}.

One consequence of giving a NCSR a vector bundle structure is topology
change. The idea that quantum gravity should lead to topology change
is not new, but up to now has been quite vague.  The topology change
described here is quite simple. It consists of two or more disjoint
surfaces each topologically equivalent to the sphere coalescing to
form one such surface. Looking at figure \ref{fig_top} we see that as
a curve is raise and lowered the corresponding surface of rotation has
a different number of disjoint connected components. This is an
example from Morse theory.

%%%%%%%%%%%%%%%%%%%%%%%%%%%%%%%%%%%%%%%%%%%%%%%%%%%%%%%%%%%%%%%%%%%%%%

\subsection{Structure of article}

We start in section \ref{ch_rev} with a review of NCSR. It is
necessary to define three separate but related algebra $\Calg$,
$\Ralg$ and $\Nalg$. It is $\Ralg$ which is closest to the algebra
defined in \cite{Gratus7}. All three algebras are defined with respect
to a $C^1$ real function $\rho$ and a parameter $\varepsilon>0$. For
$\Ralg$ a third parameter $R\in\Real$ is given. 
In section \ref{ch_rep} we give the finite dimensional unitary
representations of these three algebras.

When $\varepsilon=0$ the algebras $\Calg$, $\Ralg$ and $\Nalg$ become
commutative algebras. In section \ref{ch_top} we give a topological
meaning to the algebra $\Ralg$. This algebra forms a dense subalgebra
of $C(\man,\Cmpx)$ the algebra of complex valued continuous functions
from the surface of rotation $\man$. The surface $\man$ which depends
on $\rho$ and $R$ is a collection of disjoint surface, each surface
topologically equivalent to the sphere.  If the curve $\rho$ has more
than one local minima then the number of disjoint surface and hence
topology of $\man$ will depend on $R$. This can be seen in figure
\ref{fig_top}.

The limit as $\varepsilon\to0$ of the commutator in
$\Ralg(\rho,\varepsilon)$ gives $\man$ a Poisson structure. In fact
this structure is symplectic. This is calculated in section
\ref{ch_poi}, and this is used to write the exterior derivative,
metric, hodge dual and Laplace equation in a form easiest to convert
to the noncommutative case.

\vspace{1 em}

Throughout sections \ref{ch_triv} and \ref{ch_vec} we assume that
$\rho$ is trivial, that is, amongst other requirements, that it has just
one local minima.  Thus all corresponding surfaces $\man(\rho,R)$ will
be connected and there is no change in topology to different $R$.
However since $\rho$ is trivial we can define the algebra $\Balg$,
which depends on $\rho$ and $\varepsilon>0$. An important subalgebra
of $\Balg$ is $\Nalg$.

The representation of $\Balg$, given in section \ref{ch_trivrep},
is a Hilbert space $\Gspace$ called a trivial multi-topology
lattice, since it may be thought of as a two dimensional lattice.
This lattice is decomposed into the direct some of vector spaces
$V_n$. Each $V_n$ is an $n$ dimensional unitary representation of the
algebras $\Nalg$ and $\Calg$. Also $V_n$ is a unitary representation
of $\Ralg$  if $R$ has a certain value. 

In section \ref{ch_sym} we show how, if $\rho(z)=z^2$ then the
algebra $\Ralg$ is isomorphic to the algebra $su(2)$, whilst the
algebra $\Balg$ is equivalent to the the product of two
Heisenberg-Weyl algebras. Thus the embedding of $\Ralg\subset\Balg$
corresponds to the Jordan-Schwinger representation of $su(2)$.  In
\cite{Gratus6} we have used this Jordan-Schwinger representation to
construct a space of vector modules over $\Ralg(\rho,\varepsilon)$ for
$\rho=z^2$.  In section \ref{ch_vec} this process is repeated for a
general trivial $\rho$. We first define the
nonassociative algebra $(\Palg,\mu)$ which because it is
nonassociative the product $\mu$ is given explicitly. This is
decomposed into a set of right modules $\Pmod{r}$ over $\Ralg$, where
$\Pmod{0}=\Ralg$, and $r\in\Intg$.  Using the results of section
\ref{ch_poi}, in section \ref{ch_onef} we show how to interpret
$\Pmod{-2}$ as the space of 1-forms over $\man$ and construct an
exterior derivative $d$, which satisfies Leibniz only in the
commutative limit.  We also interpret $\Pmod{2}$ as the analogous
space (module) of tangent vector fields $T\man$. This defines a vector
as an object which can be multiplied by a scalar to give a vector, but
which is a derivative only in the commutative limit.  In section
\ref{ch_inn} we show that any operator that obeys Leibniz must be
inner, and thus the fact the our operator does not obey Leibniz must
be accepted.

\vspace{1 em}

In section \ref{ch_notriv} we return our attention to the more general
$\rho$. We can no longer form the algebra $\Balg$, however we can
still investigate the multi-topology lattices $\Gspace$. Here, once
again, $\Gspace$ is a direct sum $\oplus_s V_s$ where each $V_s$ is a
representation of $\Nalg$. Since $\Gspace$ encodes the different
topologies of $\man$, this justifies the name multi-topology lattice.

Unfortunately problems may occur near the topology change, and this
requires one of four compromises to be made. The detail of the
construction of $\Gspace$ and the compromises is given in sections
\ref{ch_G0} to \ref{ch_rep_topch}, where we demonstrate that one can
always construct a multi-topology lattice.  In section \ref{ch_topch}
we describe some operators which can be interpreted as operators for
topology change. We indicate that a topology change can be thought of
as a block diagonal matrix.  Finally in section \ref{ch_GMTL} we
propose an operator algebra which may model the dynamics of surfaces
which, although they remain axially symmetric, change shape and
interact.

\vspace{1 em}

Finally in section \ref{ch_disc} we discuss some of the areas of
research that follow from this article.

%%%%%%%%%%%%%%%%%%%%%%%%%%%%%%%%%%%%%%%%%%%%%%%%%%%%%%%%%%%%%%%%%%%%%%
%%%%%%%%%%%%%%%%%%%%%%%%%%%%%%%%%%%%%%%%%%%%%%%%%%%%%%%%%%%%%%%%%%%%%%
%%%%%%%%%%%%%%%%%%%%%%%%%%%%%%%%%%%%%%%%%%%%%%%%%%%%%%%%%%%%%%%%%%%%%%

\section{Review of Noncommutative Surfaces of Rotation (NCSR)}
\label{ch_rev}

For the purposes of this article we will review three closely related
but different algebra $\Calg,\Ralg,\Nalg$.  All three can be referred
to as the NCSR so we will use the correct symbol if we need to be
precise.  The algebra $\Ralg$ is equivalent the algebra given in
\cite{Gratus7} where NCSR were first introduced.

Let us define the domains $\Xdom=C^1(\Real\mapsto\Cmpx)$ and
$\NXdom=C^1(\Real\times\Real\mapsto\Cmpx)$. We have chosen $C^1$ but
similar results exist for $C^k$ or $C^\omega$.  We note that both
these domains are commutative algebras where the product $fg$ is the
pointwise multiplication.

The algebras $\Calg,\Ralg,\Nalg$ are all defined with respect to a
function $\rho$. Throughout this article we will assume that 
%[
\begin{align}
\begin{minipage}{10 cm}
\begin{jglist}
\item $\rho\in C^1(\Real\mapsto\Real)$, 
\item $\rho\to\infty$ as $x\to\pm\infty$,
\item $\rho$ has a finite number of turning points
\item there is no interval in which $\rho$ is a constant.  
\end{jglist}
\end{minipage}
\label{rev_rho_asum}
\end{align}
%]
All three algebras also require that we specify
$\varepsilon\in\Real$ with $\varepsilon\ge 0$. There are some
results which can be reformulated for negative or even complex
$\varepsilon$ but for this article we will assume that
$\varepsilon\ge0$. For the algebra $\Ralg$ there is a third parameter
$R\in\Real$ for which we make no further assumptions.

%%%%%%%%%%%%%%%%%%%%%%%%%%%%%%%%%%%%%%

The three algebras are given in table \ref{rev_tbl} together with a
fourth algebra $\NCalg$ which is used to relate the other three
algebras to each other. We shall call the elements $X_+$ and $X_-$ the
ladder operators, and the elements $X_0$ and $N_0$ the diagonal
operators. These names come from representations.
The expression $X_\pm^r$ means
%[
\begin{align}
X_\pm^r &= 
\begin{cases}
(X_+)^r & \text{ if } r\ge 0\\
(X_-)^{-r} & \text{ if } r< 0
\end{cases}
\label{rev_def_Xpm}
\end{align}
%]
We see from the list of generators for each algebra that the ladder
operators $X_+,X_-$ are handle differently from the diagonal operators
$X_0,N_0$. In general any $\Xdom$ function of $X_0$ or $\NXdom$
function of $(X_0,N_0)$ are allowed, but only polynomials of $X_+$ and
$X_-$. We have to handle $\Xdom$ function of $X_0$ because
$\rho\in\Xdom$ and $\rho$ appears in the defining equations of the
algebra. If we were to allow $\Xdom$ functions of $X_+$ for example
this would cause problems with the existence or otherwise of limits.

\begin{lemma}
For all four algebras 
%[
\begin{align}
[X_0,X_\pm] = \pm\varepsilon X_\pm
\label{rev_com_X0Xpm}
\end{align}
%]
whilst for $\Nalg$ and $\NCalg$
%[
\begin{align}
[N_0,X_\pm] = 0
\label{rev_com_N0Xpm}
\end{align}
%]
for $\Ralg$ we have
%[
\begin{align}
X_-X_+ &= \rho(R)-\rho(X_0+\varepsilon)
\label{rev_Ra_XmXp}
\end{align}
%]
and for $\Nalg$ we have
%[
\begin{align}
X_-X_+ &= \rho(N_0)-\rho(X_0+\varepsilon)
\label{rev_Na_XmXp}
\end{align}
%]
Every element can be written uniquely as (\ref{rev_Ca_ge}),
(\ref{rev_Ra_ge}), (\ref{rev_Na_ge}), or (\ref{rev_NCa_ge}).
This form is known as normal ordering. 
\end{lemma}

\begin{proof}
(\ref{rev_com_X0Xpm}) and (\ref{rev_com_N0Xpm}) follow from the
respective quotient equations.  (\ref{rev_Ra_XmXp}) and
(\ref{rev_Na_XmXp}) follow by considering $X_+X_-X_+$.
\end{proof}

The algebra $\NCalg$ may be thought of as the extensions of $\Calg$
with the central element $N_0$. This algebra is defined so we have the
follow lemma.

\begin{lemma}
\label{lm_rev_com_diag}
The relationship between the three algebras for a NCSR is given by the
following diagram:
%[
\begin{align}
\setlength{\unitlength}{8 mm}
\begin{picture}(7.5,2.5)
\put(0,2){\makebox(0,0){$\Calg$}}
\put(1,2){\makebox(0,0){$\hookrightarrow$}}
\put(2,2){\makebox(0,0){$\NCalg$}}
\put(4,2){\makebox(0,0){$\Nalg$}}
\put(4,0){\makebox(0,0){$\Ralg$}}
\put(2.5,2){\vector(1,0){1}}
\put(0,1.5){\vector(2,-1){3}}
\put(4,1.5){\vector(0,-1){1}}
\put(3,2.2){\makebox(0,0){${}^{q_1}$}}
\put(1.5,1){\makebox(0,0){${}^{q_2}$}}
\put(4.3,1){\makebox(0,0){${}^{q_3}$}}
\end{picture}
\label{rev_com_diag}
\end{align}
%]
where the hooked arrows refer to the natural embedding, 
$q_1$ is the quotient $X_+X_- - \rho(N_0)-\rho(X_0)\sim 0$, 
$q_2$ is the quotient $X_+X_- - \rho(R)-\rho(X_0)\sim 0$, and
$q_3$ is the quotient $N_0-R\sim 0$.
\end{lemma}

\begin{proof}
Trivial. 
\end{proof}

%%%%%%%%%%%%%%%%%%%%%%%%%%%%%%%%%%%%%%%%%%%%%%%%%%%%%%%%%%%%%%%%%%%%%%

\subsection{``Standard'' Representations of $\rho$}
\label{ch_rep}

A {\bf standard unitary representation} of $\Calg$, $\Ralg$ is $\Nalg$
defined with respect to the pair $(J,V)$ where $J=\{z\ |\ z\smin{}\le
z \le z\smax{}\}\subset\Real$ is an interval, called the {\bf
representation interval} such that
%[
\begin{equation}
\begin{aligned}
{}&|J|/\varepsilon = (z\smax{}-z\smin{})/\varepsilon\in\natnum 
\\
{}&\rho(z\smin{})=\rho(z\smax{}) 
\\
{}&\rho(z) \le \rho(z\smax{})\ \forall z\in J
\end{aligned}
\label{rep_def_J}
\end{equation}
%]
and $V$ is a finite dimensional vector space with dimension
$\dim(V)=|J|/\varepsilon$. The basis of $V$ is $\Rvec(m)$ where
$m\in\Intg$ and $m\smin{}\le m\le m\smax{}$. Here $m\smin{}$ is an
arbitrary integer and $m\smax{}=m\smin{}+\dim V -1$.

The standard unitary representation of 
$\Calg$ with respect to the pair $(J,V)$ is given by
%[
\begin{equation}
\begin{aligned}
f(X_0)\Rvec(m) &= 
f(z\smin{} - \varepsilon m\smin{} + \varepsilon m) \Rvec(m) 
\\
X_+\Rvec(m) &= 
(\rho(z\smin{})-
 \rho(z\smin{} - \varepsilon m\smin{} + \varepsilon m + \varepsilon)
)^\scrhalf 
\Rvec(m+1) 
\\
X_-\Rvec(m) &= 
(\rho(z\smin{})-
\rho(z\smin{} - \varepsilon m\smin{} + \varepsilon m))^\scrhalf 
\Rvec(m-1) 
\end{aligned}
\label{rep_stand_rep}
\end{equation}
%]

The standard representation of $\Nalg$ and $\NCalg$ with respect to
$(J,V)$ is given by 
%[
\begin{align}
f(X_0,N_0)\Rvec(m)=
f(z\smin{} - \varepsilon m\smin{} + \varepsilon m,z\smin{})\Rvec(m) 
\label{rep_rep_N0}
\end{align}
%]
and the last two equations of (\ref{rep_stand_rep}).

There exists a standard representation of $\Ralg(\rho,\varepsilon,R)$
with respect to $(J,V)$ if and only if $R=z\smin{}$. In this case the
representation is given again by (\ref{rep_stand_rep}).

Therefore unlike $\Calg$ and $\Nalg$ where there exists an infinite
set of standard representation, for $\Ralg$ there exists either one or
zero standard representation. The infinite set of standard
representation of $\Nalg$ is used in the constructions of the
multi-topology lattice representation of the NCSR. See section
\ref{ch_trivrep} for the representation of trivial NCSR, and section
\ref{ch_notriv} for the representation of non-trivial NCSR.

%%%%%%%%%%%%%%%%%%%%%%%%%%%%%%%%%%%%%%%%%%%%%%%%%%%%%%%%%%%%%%%%%%%%%%

\subsection{Topology of $\rho$}
\label{ch_top}

For a given $R\in\Real$, let $\man=\man(\rho,R)$ be the surface in
$\Real^3$
%[
\begin{align}
\man=\man(\rho,R) = \{(x^1,x^2,x^3)\in\Real^3\ |\
(x^1)^2+(x^2)^2=\rho(R)-\rho(x^3) \}
\label{top_imbed}
\end{align}
%]
If $\rho(R)$ is the value of a maxima of $\rho$ then the set of points
$\man$ obeying (\ref{top_imbed}) do not form a manifold, since two of
the surfaces are glued at a point.  If $\rho(R)$ is the value of a
minima of $\rho$ then $\man$ describes a manifold but one of the
disjoint pieces is simply a point.  The set of $R$ where $\rho(R)$ is
the value of a maxima or minima of $\rho$ are called singular. Since
we are interested in describing surfaces we must assume that $R$ is
nonsingular.

We give a coordinate system for $\man$ as $(\phi,z)$ where
%[
\begin{align}
z\in\{z\ |\ \rho(R)\ge \rho(z)\}
\label{top_z_range}
\end{align}
%]
and $0\le \phi < 2\pi$ and give the coordinate chart
$:(\phi,z)\mapsto\man$ as
%[
\begin{align}
x^1 &= (\rho(R)-\rho(z))^\scrhalf \cos\phi &
x^2 &= (\rho(R)-\rho(z))^\scrhalf \sin\phi &
x^3 &= z
\label{top_man}
\end{align}
%]
For a given nonsingular $R$, the surface $\man(\rho,R)$ is the
disjoint union of connected surfaces, each one topologically
equivalent to the sphere. The surfaces are in one to one
correspondence with the intervals in the set (\ref{top_z_range}).

We have followed the standard procedure of noncommutative geometry by
saying the algebra $\Ralg(\rho,0,R)$ is ``equivalent'' to the
manifold $\man(\rho,R)$, and that $\Ralg(\rho,\varepsilon,R)$ for
$\varepsilon\ne0$ is the noncommutative analogue of $\man(\rho,R)$ or,
alternatively, that $\Ralg(\rho,\varepsilon,R)$ is an $\varepsilon$
perturbation of $\man(\rho,R)$.  Of course an algebra is not equivalent
to a manifold. What we mean is that there exists the map
%[
\begin{align}
:\Ralg &\mapsto  C^0(\man\mapsto\Cmpx) 
\notag
\\
\sum_r X_\pm^r f_r(X_0)
&\mapsto
\sum_r e^{i r\phi} f_r(z)
\label{top_def_Xpm}
\end{align}
%]
The set $C^0(\man\mapsto\Cmpx)$ is given a commutative algebraic
structure via pointwise multiplication, so that (\ref{top_def_Xpm}) is
a homomorphism. (\ref{top_def_Xpm}) is also and injection, so we can
view $\Ralg(\rho,0,R)\subset C^0(\man\mapsto\Cmpx)$.  We give
$C^0(\man\mapsto\Cmpx)$ the $L^\infty$ topology so $\Ralg$ forms a
dense subset of $C^0(\man\mapsto\Cmpx)$.  Thus we can approximate any
continuous function on $\man$ with a series of functions in $\Ralg$
and then deform these functions to give elements in
$\Ralg(\rho,\varepsilon,R)$.

\vspace{1 em}

We define the {\bf Topology intervals} $\topI{t}$ which is a set of
intervals corresponding to a single connected surface as $\rho(R)$
moves up and down. The topology intervals are labelled $\topI{t}$
where $t\in\T$ and $\T$ is a finite index set.  These are defined as:
\begin{jglist}
\item
If $I=\{z|z\smin{}\le z\le z\smax{}\}$ then $I\in\topI{t}$ for some
$t\in\T$ if and only if 
%[
\begin{align}
\rho(z\smin{})&=\rho(z\smax{}) 
\text{\qquad and \qquad}
\rho(z)<\rho(z\smin{}) \qquad\forall\, z\in I\backslash\{z\smin{},z\smax{}\}
\label{top_def_I}
\end{align}
%]

\item
If $I\in\topI{t}$ and $I'\in\topI{t}$ then either 
$I\subset I'$ or $I'\subset I$ and there are no maxima in the sets
$I\backslash I'$ or $I'\backslash I$.
\end{jglist}

In figure \ref{fig_top} the curve $\rho$ has 5 topology intervals,
with $I_1$, $I_2$, $I_3$ and $I_4$ belonging to different
intervals. $I_3$ and $I_5$ belong to the same topology interval.

\begin{lemma}
The topology interval $\topI{t}$ may be described uniquely with
respect to four parameters
$z^{(t)}_1$, $z^{(t)}_2$, $z^{(t)}_3$, $z^{(t)}_4\in\Real\union{\pm\infty}$
which satisfy the following
\begin{jglist}
\item
$z^{(t)}_1\le z^{(t)}_2\le z^{(t)}_3\le z^{(t)}_4$, 

\item
$\rho(z^{(t)}_1)=\rho(z^{(t)}_4)$, and $\rho(z^{(t)}_2)=\rho(z^{(t)}_3)$, 

\item
$\rho$ is decreasing between $z^{(t)}_1$ and $z^{(t)}_2$

\item
$\rho$ is increasing between $z^{(t)}_3$ and $z^{(t)}_4$

\item
$\rho(z)\le\rho(z^{(t)}_2)$ for all $z^{(t)}_2\le z\le z^{(t)}_3$

\item Either $-z^{(t)}_1=z^{(t)}_4=\infty$ or $z^{(t)}_1$ is a local
maxima or $z^{(t)}_4$ is a local maxima

\item Either $z^{(t)}_2=z^{(t)}_3$ or there exists a point
$z^{(t)}_3<z^{(t)}_5<z^{(t)}_4$ such that
$\rho(z^{(t)}_5)=\rho(z^{(t)}_3)$
\end{jglist}

$\topI{t}$ is now defined as the set of intervals
%[
\begin{align}
&\topI{t}=\left\{ I=\{z|z\smin{}\le z\le z\smax{}\} \ \bigg|\
\rho(z\smin{})=\rho(z\smax{})\,,\
z^{(t)}_1\le z\smin{}<z^{(t)}_2\,,\
z^{(t)}_3<z\smax{} \le z^{(t)}_4  \right\}
\end{align}
%]
\end{lemma}

\begin{proof}
Trivial. (See Figure \ref{fig_topi}.)
\end{proof}

\begin{lemma}
The number of topology intervals is bounded by
%[
\begin{align}
|\T| \le
\text{number of maxima + number of minima}
\end{align}
%]
with equality when all the maxima have different values.
\end{lemma}

\begin{proof}
If the number of maxima is finite then there is one topology interval
which extend to infinity. Each maxima now creates two new topology
intervals, hence result. 
\end{proof}

As $R$ increases and crosses a singular point, the topology of the
manifold undergoes a transition, as two or more adjacent intervals
coalesce or a single interval splits.  This corresponds to two or more
surfaces coalescing or a single surface bifurcating.  These changes
can be encoded into a function $\piT:\T\mapsto\T\union\{\infty\}$
where $\topI{\piT(t)}$ is directly above $\topI{t}$ or
$\topI{t}=\infty$ if it is the highest topology interval. 
Thus 
%[
\begin{align}
\rho(z^{(t_2)}_2) = \rho(z^{(t_1)}_1)
\qquad
\text{when}
\quad
t_2 = \piT(t_1)
\end{align}
%]
If $\piT^{-1}\{s\}=\{t_1,t_2,\ldots,t_N\}$ then there is a
transition where the topology intervals
$\topI{{t_1}},\topI{{t_2}},\ldots,\topI{{t_N}}$ coalescing into the
interval $\topI{s}$. This implies there are $N-1$ maxima all of the
same value, one maxima between successive $\topI{t_i}$.

For example in figure \ref{fig_top}, if we let $I_1\in\topI{1}$,
$I_2\in\topI{2}$, $I_3,I_5\in\topI{3}$, $I_4\in\topI{4}$ and let
$\topI{5}$ be the unbounded topology interval then
$\piT(1)=\piT(2)=4$, $\piT(3)=\piT(4)=5$ and
$\piT(5)=\infty$.

%%%%%%%%%%%%%%%%%%%%%%%%%%%%%%%%%%%%%%%%%%%%%%%%%%%%%%%%%%%%%%%%%%%%%%

\subsection{The metric and differential structure of $\man$
in terms of the Poisson structure}
\label{ch_poi}

The limit of the noncommutative structure gives $\man$ a Poisson
structure \cite{Gratus7}:
%[
\begin{align}
\PB(f,h) &= 
\lim_{\varepsilon\to0}
\left( \frac1\varepsilon
[f,h]
\right) 
=
-i\left(
\frac{\partial f}{\partial z}\frac{\partial h}{\partial \phi}
-
\frac{\partial h}{\partial z}\frac{\partial f}{\partial \phi}
\right)
\label{poi_def}
\end{align}
%]

We can write some of the standard objects of differential geometry
simply in terms of the Poisson structure, the elements of
$C(\man,\Cmpx)$ and the elements $dX_+$, $dX_-$ and $dX_0$. This is 
useful since these expression are the easiest to extend to the
noncommutative case. Although in this article we will attempt only a
definition of $df$ by extending (\ref{poi_df}).
%[
\begin{align}
df 
&= 
\frac{1}{\rho(R)-\rho(X_0)}
X_-\bigg(dX_+\PB(X_0,f) - dX_0\PB(X_+,f)\bigg)
\label{poi_df}
\\
df
&=
\frac{-1}{\rho(R)-\rho(X_0)}
X_+\bigg(dX_-\PB(X_0,f) - dX_0\PB(X_-,f)\bigg)
\label{poi_df2}
\end{align}
%]

The metric $g$ on $\man$ is given by the pull back of the
metric $\Real^3$ using the mapping given by (\ref{top_man}). 
%[
\begin{align}
g &=
dz\otimes dz + \tfrac12 dX_+\otimes dX_- + \tfrac12 dX_-\otimes dX_+ 
\notag\\
&= 
\frac{C(z)}{4(\rho(R)-\rho(z))}
dz\otimes dz +
(\rho(R)-\rho(z)) d\phi\otimes d\phi
\label{poi_g}
\end{align}
%]
where
%[
\begin{align}
C(z)=\rho'(z)^2+4\rho(R)-4\rho(z)
\label{poi_C}
\end{align}
%]
Using $g$ we construct the map $\wt{}:T^*\man\mapsto T\man$. Thus we
can express the metric $g$ solely in terms of functions on $\man$,
i.e. elements of $C(\man,\Cmpx)$.
%[
\begin{align}
g(\wt{df},\wt{dh}) &=
\star^{-1}(df\wedge\star dh)
=
\frac{-2}{C(X_0)}
\bigg(
\PB(X_+,f)\PB(X_-,h) +
\PB(X_-,f)\PB(X_+,h) +
2\PB(X_0,f)\PB(X_0,h) 
\bigg)
\label{poi_gfh}
\end{align}
%]
The metric $g$ also defines a Hodge dual
$\star:\Lambda(\man)\mapsto\Lambda(\man)$. This gives
%[
\begin{align}
\star df &=
\frac{2i}{C(X_0)^\scrhalf}
\bigg(
dX_-\PB(f,X_+) + dX_+\PB(f,X_-) + 2dX_0\PB(f,X_0)
\bigg)
\label{poi_sdf}
\end{align}
%]
and the Laplace operator
%[
\begin{align}
\star^{-1}(d\star df) &= 
\frac{-2}{C(X_0)}
\bigg(\PB(X_-,{\PB(X_+,f)})+\PB(X_+,{\PB(X_-,f)})
+2\PB(X_0,{\PB(X_0,f)}) \bigg) 
\notag\\
&+
\frac{2}{C(X_0)^2}
\rho'(X_0)(\rho''(X_0)-2)
\bigg(X_-\PB(X_+,f)-X_+\PB(X_-,f)\bigg)
\label{poi_lap}
\end{align}
%]

As stated in the introduction we wish to find corresponding
definitions for these geometric objects, but for $f\in\Ralg$ instead
of $f\in C^0(\man,\Cmpx)$.  The idea is that these new definitions
reduce to the above expressions when $\varepsilon\to0$.  In this
article, we only suggest a new definitions to the exterior derivative,
which reduces to (\ref{poi_df}) and (\ref{poi_df2}) in the limit
$\varepsilon\to0$.

%%%%%%%%%%%%%%%%%%%%%%%%%%%%%%%%%%%%%%%%%%%%%%%%%%%%%%%%%%%%%%%%%%%%%%
%%%%%%%%%%%%%%%%%%%%%%%%%%%%%%%%%%%%%%%%%%%%%%%%%%%%%%%%%%%%%%%%%%%%%%
%%%%%%%%%%%%%%%%%%%%%%%%%%%%%%%%%%%%%%%%%%%%%%%%%%%%%%%%%%%%%%%%%%%%%%
%%%%%%%%%%%%%%%%%%%%%%%%%%%%%%%%%%%%%%%%%%%%%%%%%%%%%%%%%%%%%%%%%%%%%%

\section{Trivial NCSR: The algebra $\Balg$}
\label{ch_triv}

Throughout this section we will assume that $\rho$ is trivial.
This means that as well as the assumptions given by
(\ref{rev_rho_asum}) we also assume
%[
\begin{align}
\begin{minipage}{10 cm}
\begin{jglist}
\item There exists a single minima at $z_0$.
\item $\rho''(z_0)\ne 0$
\item $\rho'(z)\ne0$ for all $z\ne z_0$
\end{jglist}
\end{minipage}
\label{triv_rho_asum}
\end{align}
%]
Such a $\rho$ will have just one topology interval and in the
classical limit $\man(\rho,R)$, for all $R\ne0$, is connected surface
topologically equivalent to the sphere. Thus there is no topology
change associated with a $\rho$ obeying (\ref{triv_rho_asum}).

We set $D=\rho(z_0)$ and define the functions
%[
\begin{align}
&\tau:\Real\mapsto\Real \qquad
\tau(x)=(\rho(x)-D)^\scrhalf \qquad
\text{and $\tau$ is decreasing}
\label{triv_def_tau}
\\
&\omega:\Real\mapsto\Real \qquad
\text{implicitly by} \qquad
\tau(\omega(x))+\tau(\omega(x)+ x) = 0
\label{triv_def_om}
\end{align}
%]
For $k\in\Intg$ we define the operators $N_k$ and $X_k$ 
(where $X_1$ and $X_{-1}$ are not to be confused with $X_+$ and $X_-$)
via
%[
\begin{align}
N_{k} = \omega(\omega^{-1}(N_0)+\varepsilon k) 
\qquad\text{and \qquad}
X_{k} = N_k - N_0 + X_0
\label{triv_def_Nk}
\end{align}
%]
\begin{lemma}
The functions $\tau$ and $\omega$ and their inverses $\tau^{-1}$ and
$\omega^{-1}$ are all well defined, strictly decreasing and belong to
$\Xdom$.
The operator $N_k$ is an $\Xdom$ function of $N_0$ whilst the
operator $X_k$ is a $\NXdom$ function of $(X_0,N_0)$.
\end{lemma}

\begin{proof}
Follows from showing of $\tau'(z)<0$ and $\omega'(z)<0$ for all
$z\in\Real$.
\end{proof}

Let us define the algebra $\Balg=\Balg(\rho,\varepsilon)$ as that
generated by $a_+,a_-,b_+,b_-$ and $\NXdom$ functions of $(X_0,N_0)$ 
quotiented by the following relationships:
%[
\begin{align}
[X_0,N_0]&=0 
& 
f(X_{j},N_{k}) a_\pm &=
a_\pm f(X_{j\pm1}\pm\varepsilon,N_{k\pm1}) 
&
f(X_{j},N_{k}) b_\pm &=
b_\pm f(X_{j\pm1},N_{k\pm1}) 
\label{triv_com_NX}
\end{align}
%]
%[
\begin{equation}
\begin{aligned}
a_-a_+ &= \tau(N_0)-\tau(X_0+\varepsilon)
&
b_-b_+ &= \tau(N_0)+\tau(X_0)
\\
a_+a_- &= \tau(N_{-1})-\tau(X_{-1})\qquad
&
b_+b_- &= \tau(N_{-1})+\tau(X_{-1})
\end{aligned}
\label{triv_prod_aa}
\end{equation}
%]
%[
\begin{equation}
\begin{aligned}
b_+a_+ &= a_+b_+ 
\left( \frac{
(\tau(N_0)-\tau(X_0+\varepsilon))
(\tau(N_1)+\tau(X_1+\varepsilon))
}{
(\tau(N_1)-\tau(X_1+\varepsilon))
(\tau(N_0)+\tau(X_0))
}\right)^\scrhalf
\\
b_-a_- &= a_-b_- 
\left( \frac{
(\tau(N_{-1})-\tau(X_{-1}))
(\tau(N_{-2})+\tau(X_{-2}-\varepsilon))
}{
(\tau(N_{-2})-\tau(X_{-2}))
(\tau(N_{-1})+\tau(X_{-1}))
}\right)^\scrhalf
\\
b_-a_+ &= a_+b_- 
\left( \frac{
\rho(N_0)-\rho(X_0+\varepsilon)
}{
(\tau(N_{-1})-\tau(X_{-1}+\varepsilon))
(\tau(N_{-1})+\tau(X_{-1}))
}\right)^\scrhalf
\\
b_+a_- &= a_-b_+ 
\left( \frac{
(\tau(N_{-1})-\tau(X_{-1}))
(\tau(N_{-1})+\tau(X_{-1}-\varepsilon))
}{
\rho(N_0)-\rho(X_0)
}\right)^\scrhalf
\end{aligned}
\label{triv_com_ab}
\end{equation}
%]

%%%%%%%%%%%%%%%%%%%%

The elements $a_+,a_-,b_+,b_-$ are known as hopping operators. This
name comes from the two dimensional lattice representation that we
discuss in the following section.  There is a Hermitian conjugate on
$\Balg$ given by
%[
\begin{align}
&\dagger:\Balg\mapsto\Balg \qquad
(a_\pm)^\dagger=a_\mp \qquad
(b_\pm)^\dagger=b_\mp \qquad
(f(X_0,N_0))^\dagger=\cnj{f}(X_0,N_0)
\quad\text{where}\quad
\cnj{f}(x,y)=\cnj{f(\cnj{x},\cnj{y})}
\notag
\\
&
\text{\qquad and \qquad}
(\xi\zeta\lambda)^\dagger=\cnj{\lambda}\zeta^\dagger \xi^\dagger\,,\ 
\text{\qquad for \qquad}
\lambda\in\Cmpx\,,\ \xi,\zeta\in\Balg
\label{triv_dag}
\end{align}
%]
The algebra $\Balg$ satisfies the
following theorem:

\begin{theorem}
\label{th_gen_B}
The definition of $\Balg$ is consistent with (\ref{triv_dag}).  All
the subexpression in (\ref{triv_com_ab}) which don't contain
$a_\pm,b_\pm$ are $\NXdom$ functions of $(X_0,N_0)$ and are real and
strictly positive. Thus the commutation relations are well defined.
The general element of $\Balg$ can be written
%[
\begin{align}
\xi = \sum_{rs}a_\pm^r b_\pm^s \xi_{rs}(X_0,N_0)
\label{triv_B_gen_el}
\end{align}
%]
where the sum is finite, $\xi_{rs}\in\NXdom$ and $a_\pm^r$ and $b_\pm^r$
are defined as in (\ref{rev_def_Xpm}).

The algebra $\Nalg$ is a subalgebra of $\Balg$ where we set 
%[
\begin{align}
X_+=b_-a_+ 
\text{\qquad and \qquad} 
X_-=a_-b_+ 
\label{triv_def_Xpm}
\end{align}
%]
Given $\xi\in\Balg$, the following are equivalent
%[
\begin{align}
&\bullet \qquad \xi \in\Nalg 
\\
&\bullet \qquad [N_0,\xi]=0 
\label{triv_N_com}
\\
&\bullet \qquad \text{we can write $\xi$ as }
\sum_{r}a_\pm^r b_\pm^{-r} \xi_{r}(X_0,N_0)
\label{triv_N_sum}
\end{align}
%]
\end{theorem}

\begin{proof}
We show that (\ref{triv_dag}) is constant with
(\ref{triv_com_NX}-\ref{triv_com_ab}) by direct substitution. We see
$(N_k)^\dagger=N_k$, so for example
%[
\begin{align*}
(N_k a_\pm)^\dagger =
a_\mp N_k =
N_{k\pm 1} a_\mp =
(a_\pm N_{k\pm 1} )^\dagger 
\end{align*}
%]

To show that the expressions in $(X_0,N_0)$ on the
right hand side of (\ref{triv_com_ab}) are in $\NXdom$ and strictly
positive consider the subexpression
%[
\begin{align}
f(X_0,N_0) &= \frac{\tau(N_1)+\tau(X_1+\varepsilon)}
{\tau(N_0)+\tau(X_0)}
\notag
\end{align}
%]
which occurs in the first equation of (\ref{triv_com_ab}), this
may be expressed
%[
\begin{align}
f(z,y) &= \frac{\tau(y^+)+\tau(y^+-y+z+\varepsilon)}
{\tau(y)+\tau(z)}
\notag
\end{align}
%]
where $z=X_0$, $y=N_0$ and
$y^+=\omega(\omega^{-1}(y)+\varepsilon)$. Using (\ref{triv_def_om}) we
see
%[
\begin{align}
f(z,y) &= \frac{\tau(y^+-y+z+\varepsilon)-
\tau(y^+ + \omega^{-1}(y)+\varepsilon)}
{\tau(z)-\tau(y+\omega^{-1}(y))}
\notag
\end{align}
%]
Since $\tau$ is single valued then the numerator and denominator of
the above equation are zero when $z=y+\omega^{-1}(y)$. The
value of $f$ at this point is
%[
\begin{align}
f(y+\omega^{-1}(y),y) &=
\frac{\tau'(y^++\omega^{-1}(y)+\varepsilon)}
{\tau'(y+\omega^{-1}(y))}
>0
\notag
\end{align}
%]
This result is the similar for all the fractions, hence result.

Equations (\ref{triv_prod_aa}) to (\ref{triv_com_ab}) that
define $\Balg$ are sufficient to reduce any expression in the hopping
operators and $\NXdom$ functions of $(X_0,N_0)$ into
(\ref{triv_B_gen_el}). To see this take word constructed from
generators. Use the commutation relations (\ref{triv_com_NX}) and
(\ref{triv_com_ab}) to push all the $a_+$ and $a_-$ to the left. Now
use (\ref{triv_prod_aa}) to remove all $a_+a_-$ pairs. Keep the
resulting term in $a_\pm^r$, and push the $\NXdom$ function in
$(X_0,N_0)$ to the right. Now do the same with the $b_+$ and $b_-$
terms.

\vspace{1 em}

To show that $\Nalg$ is a subalgebra of $\Balg$ we reproduce
the defining equations of
$\Nalg$. We already have $[N_0,X_0]=0$. For $f\in\NXdom$
%[
\begin{align*}
f(X_0,N_0) X_+
=
f(X_0,N_0) b_- a_+
=
b_- f(X_{-1},N_{-1}) a_+
=
b_- a_+ f(X_0+\varepsilon,N_0)
=
X_+ f(X_0+\varepsilon,N_0)
\end{align*}
%]
and likewise for $X_-$, hence (\ref{rev_Na_f}).
%[
\begin{align*}
X_+X_- &= b_-a_+a_-b_+ \\
&= b_-(\tau(N_{-1})-\tau(X_{-1})) b_+ \\
&= b_-b_+ (\tau(N_{0})-\tau(X_0)) \\
&= (\tau(N_0)+\tau(X_0))(\tau(N_0)-\tau(X_0)) \\
&= \tau(N_0)^2-\tau(X_0)^2 \\
&= (\rho(N_0)-D)-(\rho(X_0)-D) = \rho(N_0) - \rho(X_0)
\end{align*}
%]
Hence (\ref{rev_Na_XX})
Thus the subalgebra generated $X_+,X_-$ given by (\ref{triv_def_Xpm})
and $\NXdom$ function of $(X_0,N_0)$ is indeed $\Nalg$.

\vspace{1 em}

Clearly if $f\in\Nalg$ then $[N_0,f]=0$. If $[N_0,f]=0$ then write
$f$ is in  (\ref{triv_B_gen_el}). For each term in the sum,
(\ref{triv_com_NX}) implies that the number $a_+$ must equal the
number of $b_-$, and the number of $b_+$ must equal the number of 
$a_-$, thus (\ref{triv_N_sum}). 

If $f$ is written in (\ref{triv_N_sum}) then use (\ref{triv_com_ab})
to permute the $a_\pm$ and $b_\pm$ to give a sequence
$a_-b_+a_-b_+\cdots$ or $b_-a_+b_-a_+\cdots$ which is replaced with
$X_+^r$ or $X_-^r$. This gives (\ref{rev_Na_ge}).
\end{proof}

%%%%%%%%%%%%%%%%%%%%%%%%%%%%%%%%%%%%%%%%%%%%%%%%%%

If we were to admit $\Xdom$ functions of the hopping operators then
we could construct the elements $X_0$ and $N_0$ as follows:
%[
\begin{equation}
\begin{aligned}
N_0 &= 
\omega(\omega^{-1}
(\tau^{-1}(\tfrac12 b_+b_- + \tfrac12 a_+a_-))+\varepsilon) \\
X_0 &= \tau^{-1}(\tfrac12 b_+b_- - \tfrac12 a_+a_-)
-\tau^{-1}(\tfrac12 b_+b_- + \tfrac12 a_+a_-)
+\omega(\omega^{-1}
(\tau^{-1}(\tfrac12 b_+b_- + \tfrac12 a_+a_-))+\varepsilon) 
\end{aligned}
\label{triv_inv_N0X0}
\end{equation}
%]

Clearly $\varepsilon$ is a parameter for the noncommutativity since
when $\varepsilon=0$ then $\Balg(\rho,\varepsilon=0)$ reduces to a
commutative algebra. This algebra is isomorphic to a dense subalgebra
of continuous functions on $\Real^4$. To see this set $a_\pm=x^1\pm
ix^2$ and $b_\pm=y^1\pm iy^2$. This give the map
$\Balg(\rho,0)\hookrightarrow C^0(\Real^4\mapsto\Cmpx)$. If we give
$C^0(\Real^4\mapsto\Cmpx)$ the standard continuity norm then
$\Balg(\rho,0)$ forms a dense set.  The diagonal operators are given
by $N_0=\tau^{-1}(\norm{x}+\norm{y})$ and
$X_0=\tau^{-1}(\norm{x}-\norm{y})$.

The noncommutative structure of $\Balg$ gives rise to a Poisson
structure on $\Real^4$. This structure is complicated to write down in
terms of the coordinates $x^1,x^2,y^1,y^2$. 

%%%%%%%%%%%%%%%%%%%%%%%%%%%%%%%%%%%%%%%%%%%%%%%%%%%%%%%%%%%%%%%%%%%%%%

\subsection{Lattice representations of a trivial NCSR}
\label{ch_trivrep}

The reason for the complicated commutation relations given by
(\ref{triv_com_ab}) is that there is a natural representation of
$\Balg$ as a two dimensional lattice.

Let $\Gspace=\Gspace(\rho,\varepsilon)$ be a Hilbert space with
orthonormal basis $\rvec(n,m)$ with $n,m\in\Intg$, and $n\ge0$, $0\le
m\le n-1$. Let the dual of $\rvec(n,m)$ be written $\lvec(n,m)$. There
is a representation of $\Balg$ given by
%[
\begin{align}
\begin{array}{r@{}l}
f(X_0,N_0)\rvec(n,m) &= 
f(\omega(\varepsilon n)+\varepsilon m,\omega(\varepsilon n)) 
\rvec(n,m) 
\\
a_+\rvec(n,m) &= \big(
\tau(\omega(\varepsilon n)) - 
\tau(\omega(\varepsilon n) + \varepsilon m  + \varepsilon)
\big)^\scrhalf \rvec(n+1,m+1) 
\\
a_-\rvec(n,m) &= \big(
\tau(\omega(\varepsilon n-\varepsilon)) - 
\tau(\omega(\varepsilon n-\varepsilon) + \varepsilon m)
\big)^\scrhalf \rvec(n-1,m-1) 
\\
b_+\rvec(n,m) &= \big(
\tau(\omega(\varepsilon n)) + 
\tau(\omega(\varepsilon n) + \varepsilon m)
\big)^\scrhalf \rvec(n+1,m) 
\\
b_-\rvec(n,m) &= \big(
\tau(\omega(\varepsilon n-\varepsilon)) + 
\tau(\omega(\varepsilon n-\varepsilon) + \varepsilon m)
\big)^\scrhalf \rvec(n-1,m) 
\end{array}
\label{trivrep_rep}
\end{align}
%]

\begin{lemma}
The representation given above by (\ref{trivrep_rep}) is indeed a
representation and it is unitary. 
\end{lemma}
\begin{proof}
This simply involves showing all the relations
(\ref{triv_com_NX}-\ref{triv_com_ab}) are consistent with
(\ref{trivrep_rep}). This should not surprise us since the relations
where constructed so that (\ref{trivrep_rep}) was a representation.
This representation of $\Balg$ is unitary since the hermitian
conjugate of $f$ is the adjoint:
%[
\begin{align*}
\lvec(n+1,m+1)a_+\rvec(n,m)
=
\big(
\tau(\omega(\varepsilon n)) - 
\tau(\omega(\varepsilon n) + \varepsilon m  + \varepsilon)
\big)^\scrhalf
=
\lvec(n,m)a_-\rvec(n+1,m+1)
\\
\lvec(n+1,m)b_+\rvec(n,m)
=
\big(
\tau(\omega(\varepsilon n)) + 
\tau(\omega(\varepsilon n) + \varepsilon m )
\big)^\scrhalf
=
\lvec(n,m)b_-\rvec(n+1,m)
\end{align*}
%]
\end{proof}

For each $n\in\natnum$ let $V_n\subset\Gspace$ be the subspace
$V_n=\spanrm\{\rvec(n,m)\ |\ m=0,\ldots,n-1\}$ and let $J_n$ be the
interval $J_n=\{z\,|\,z\smin{n}\le z\le z\smax{n}\}$ where 
$z\smin{n}=\omega(\varepsilon n)$ and
$z\smax{n}=\omega(\varepsilon n)+ \varepsilon n$.

\begin{lemma}
The Hilbert space $\Gspace$ is also a representation of
$\Nalg\subset\Balg$ given by (\ref{trivrep_rep}) and
(\ref{triv_def_Xpm}). This is given by
%[
\begin{equation}
\begin{aligned}
f(X_0,N_0)\rvec(n,m) 
&= f(z\smin{n}+\varepsilon m,z\smin{n}) \rvec(n,m)
\\
X_+\rvec(n,m) &= 
(\rho(z\smin{n})-
 \rho(z\smin{n} + \varepsilon m + \varepsilon)
)^\scrhalf 
\rvec(n,m+1) 
\\
X_-\rvec(n,m) &= 
(\rho(z\smin{n})-
\rho(z\smin{n} + \varepsilon m))^\scrhalf 
\rvec(n,m-1) 
\end{aligned}
\label{trivrep_A}
\end{equation}
%]
For each $n\in\natnum$, the interval $J_n$ is a representation
interval, and the pair $(J_n,V_n)$ define a standard unitary
representation of $\Calg$, $\Nalg$ and
$\Ralg(\rho,\varepsilon,R=z\smin{n})$ given by (\ref{trivrep_A}).
\end{lemma}

\begin{proof}
By substituting (\ref{trivrep_rep}) into (\ref{triv_N_sum}) we get
(\ref{trivrep_A}). 

Clearly $|J_n|/\varepsilon=n\in\natnum$ and
$\rho(z\smin{n})=\rho(z\smax{n})$.  Since $\rho$ has one minima which
must lie between $z\smin{n}$ and $z\smax{n}$ then $J_n$ obeys
(\ref{rep_def_J}) and $J_n$ is a representation interval. Equations 
(\ref{trivrep_A}) coincide with (\ref{rep_stand_rep}) and
(\ref{rep_rep_N0}) when $z\smin{}=z\smin{n}$ and $m\smin{}=0$.

For the algebra $\Ralg(\rho,\varepsilon,R=z\smin{n})$ use the last
four equations in (\ref{rep_stand_rep}).
\end{proof}

\begin{figure}[t]
\epsfysize 5 cm
\centerline{\epsffile{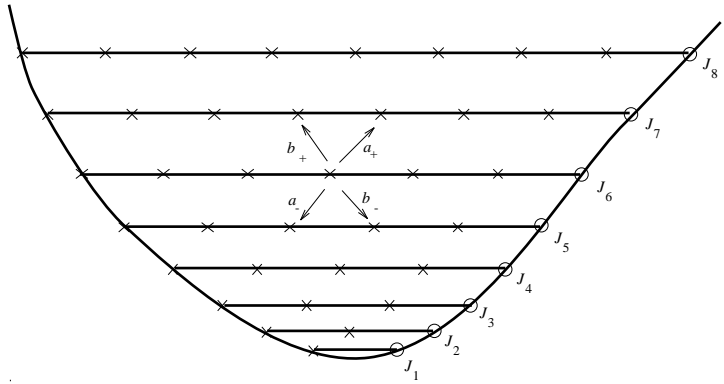}}
\caption{First few $V_n$ for a trivial $\rho$}
\label{fig_triv}
\vspace{10em}
\end{figure}

In figure \ref{fig_triv} we view each basis vector
$\rvec(n,m)$ as a point in $\Real^2$ with coordinates
%[
\begin{align}
x-\text{coord}=\lvec(n,m)X_0\rvec(n,m) 
\qquad\qquad
y-\text{coord}=\lvec(n,m)\rho(N_0)\rvec(n,m) 
\label{trivrep_coords}
\end{align}
%]
In figure \ref{fig_triv} these points are represented by
crosses. The last point on the right of each $J_n$ does not represent
a basis vector and so is drown with a circle.
Since the $y$-coordinate is independent of $m$ the $n$ points in each $V_n$
lie on the same horizontal line. This line is labelled $J_n$ although
the representation interval $J_n$ refers only to the $x$-coordinate. 

Clearly all the crosses lie on or above the curve $\rho$ since
$\rho(z\smin{n})-\rho(z)\ge 0$ for $z\in J_n$. 
The $y$-coordinate of the lines $J_n$ increases with $n$ since
%[
\begin{align}
\lvec(n+1,m)\rho(N_0)\rvec(n+1,m) \ge
\lvec(n,m)\rho(N_0)\rvec(n,m) 
\label{trivrep_inc_rhoN}
\end{align}
%]
The arrow representing $a_+$ always point to the top right whilst the
arrow representing $b_+$ always point to the top left. This is because
%[
\begin{align}
\lvec(n+1,m) X_0 \rvec(n+1,m) \le
\lvec(n,m) X_0 \rvec(n,m) \le
\lvec(n+1,m+1) X_0 \rvec(n+1,m+1) 
\label{trivrep_posn_X0}
\end{align}
%] 

%%%%%%%%%%%%%%%%%%%%

\vspace{1 em}

Although the representation is not faithful we do have the following
lemma.
\begin{lemma}
If $f(x,y,\varepsilon)$ is $\NXdom$ in $(X_0,N_0)$ and analytic in
$\varepsilon$ in a domain about $\varepsilon=0$ and if
%[
\begin{align}
f(X_0,N_0,\varepsilon)\rvec(n,m)=0
\end{align}
%]
for all $\rvec(n,m)\in\Gspace(\rho,\varepsilon)$ and for all
$\varepsilon>0$ then $f\equiv 0$.

If $\rho$ is analytic and $\xi$ is a word constructed from the
generators of $\Balg$ without explicit $\varepsilon$ and
%[
\begin{align}
\xi\rvec(n,m)=0
\end{align}
%]
for all $\rvec(n,m)\in\Gspace(\rho,\varepsilon)$ and for all
$\varepsilon>0$ then $\xi\equiv 0$.
\end{lemma}

\begin{proof}
Fix $y$. Let $y_1>y$ satisfy $\rho(y)=\rho(y_1)$. Now fix $x$ so that
$(x-y)/(y_1-x)\in\ratnum$ and let $\varepsilon_0$ be the highest
common factor of $(x-y)$ and $(y_1-x)$. So
$N=(y_1-y)/\varepsilon_0\in\natnum$ and
$M=(x-y)/\varepsilon_0\in\natnum$.

For each $t\in\natnum$ choose
$\varepsilon=\varepsilon_0/t$, $n=tN$ and $m=tM$ 
%[
\begin{align}
0=
\lvec(n,m)f(X_0,N_0,\varepsilon)\rvec(n,m)
=
f(x,y,\varepsilon)
\notag
\end{align}
%]
and since $f(x,y,\varepsilon)$ is analytic in $\varepsilon$ about
$\varepsilon=0$, this implies $f(x,y,\varepsilon)=0$ for all
$\varepsilon$. 

The conditions on $x$ implies that $f(x,y,\varepsilon)=0$ for a dense
set of $x$, and since $f$ is $\Xdom$ in $x$ it is true for all
$x$. Finally since we were free to choose $y$ we have $f\equiv0$.

If $\rho$ is analytic and $\xi$ is generated as stated then $\xi$ can
be rewritten as (\ref{triv_B_gen_el}) with the $f_{rs}\in\Xdom$ with
an implicit analytic dependence on $\varepsilon$. Thus we can apply
the first part of this lemma.
\end{proof}

%%%%%%%%%%%%%%%%%%%%%%%%%%%%%%%%%%%%%%%%%%%%%%%%%%%%%%%%%%%%%%%%%%%%%%
\subsection{Under-hopping operators}
\label{ch_trivalt}

Let us consider an alternative set of hopping operators given by
%[
\begin{equation}
\begin{aligned}
a'_+ = a_+ \left( 
\frac{\tau(N_1) - \tau(X_1+\varepsilon)}
{\tau(N_0)-\tau(X_0+\varepsilon)}
\right)^\scrhalf
\\
a'_- = a_- \left( 
\frac{\tau(N_0) - \tau(X_0)}
{\tau(N_{-1})-\tau(X_{-1})}
\right)^\scrhalf
\end{aligned}
\qquad\qquad
\begin{aligned}
b'_+ = b_+ \left( 
\frac{\tau(N_1) + \tau(X_1+\varepsilon)}
{\tau(N_{0})+\tau(X_{0})}
\right)^\scrhalf
\\
b'_- = b_- \left( 
\frac{\tau(N_0) + \tau(X_0+\varepsilon)}
{\tau(N_{-1})+\tau(X_{-1})}
\right)^\scrhalf
\end{aligned}
\label{trivalt_abp}
\end{equation}
%]
The operators $a'_+,a'_-,b'_+,b'_-$ are clearly in $\Balg$ since the
right hand side of (\ref{trivalt_abp}) are $\NXdom$ functions of
$(X_0,N_0)$. This is similar to the proof of theorem \ref{th_gen_B}.

It is easy to show that
%[
\begin{align}
X_+ = a'_+ b'_-
\text{\qquad and \qquad}
X_- = b'_+ a'_-
\end{align}
%]
Thus we could have constructed $\Balg$ using these alternative hopping
operators instead of $a_\pm,b_\pm$. All the results would have been
similar and the algebra $\Balg$ constructed using $a'_\pm,b'_\pm$ is
isomorphic to the algebra constructed using $a_\pm,b_\pm$.
This is not the case for the non-trivial NCSR as we will see in the
next section.

To distinguish between the two sets of hopping operators we will call
the set $\{a_\pm,b_\pm\}$ the over-hopping operators, whilst the set
$\{a'_\pm,b'_\pm\}$ the under-hopping operators. This is due to the
following diagram.
%[
\begin{align}
\setlength{\unitlength}{15 mm}
\begin{picture}(3,2.2)(0,-1)
\put(0,0){$\rvec(n,m)$}
\put(1,1){$\rvec(n+1,m+1)$}
\put(2,0){$\rvec(n,m+1)$}
\put(1,-1){$\rvec(n-1,m)$}
\put(0.3,0.3){\vector(1,1){0.6}}
\put(1.4,0.9){\vector(1,-1){0.6}}
\put(0.3,-0.1){\vector(1,-1){0.6}}
\put(1.4,-0.7){\vector(1,1){0.6}}
\put(0.6,0.6){\makebox(0,0)[br]{$a_+$}}
\put(1.7,0.6){\makebox(0,0)[bl]{$b_-$}}
\put(0.6,-0.4){\makebox(0,0)[tr]{$a'_+$}}
\put(1.7,-0.4){\makebox(0,0)[tl]{$b'_-$}}
\put(0.5,0.1){\vector( 1,0){1.4}}
\put(1.3,0.15){\makebox(0,0)[b]{$X_+$}}
\end{picture}
&&&
\setlength{\unitlength}{15 mm}
\begin{picture}(3,2.2)(0,-1)
\put(0,0){$\rvec(n,m)$}
\put(1,1){$\rvec(n+1,m+1)$}
\put(2,0){$\rvec(n,m+1)$}
\put(1,-1){$\rvec(n-1,m)$}
\put(0.9,0.9){\vector(-1,-1){0.6}}
\put(2.0,0.3){\vector(-1,1){0.6}}
\put(0.9,-0.7){\vector(-1,1){0.6}}
\put(2.0,-0.1){\vector(-1,-1){0.6}}
\put(0.6,0.6){\makebox(0,0)[br]{$a_-$}}
\put(1.7,0.6){\makebox(0,0)[bl]{$b_+$}}
\put(0.6,-0.4){\makebox(0,0)[tr]{$a'_-$}}
\put(1.7,-0.4){\makebox(0,0)[tl]{$b'_+$}}
\put(1.9,0.1){\vector(-1,0){1.4}}
\put(1.3,0.15){\makebox(0,0)[b]{$X_-$}}
\end{picture}
\end{align}
%]
The representation of the under-hopping operators
are given by
%[
\begin{align}
\begin{array}{r@{}l}
a'_+\rvec(n,m) &= \big(
\tau(\omega(\varepsilon n + \varepsilon)) - 
\tau(\omega(\varepsilon n + \varepsilon) + \varepsilon m  + \varepsilon)
\big)^\scrhalf \rvec(n+1,m+1) 
\\
a'_-\rvec(n,m) &= \big(
\tau(\omega(\varepsilon n)) - 
\tau(\omega(\varepsilon n) + \varepsilon m)
\big)^\scrhalf \rvec(n-1,m-1) 
\\
b'_+\rvec(n,m) &= \big(
\tau(\omega(\varepsilon n + \varepsilon)) + 
\tau(\omega(\varepsilon n + \varepsilon) + \varepsilon m + \varepsilon)
\big)^\scrhalf \rvec(n+1,m) 
\\
b'_-\rvec(n,m) &= \big(
\tau(\omega(\varepsilon n)) + 
\tau(\omega(\varepsilon n) + \varepsilon m)
\big)^\scrhalf \rvec(n-1,m) 
\end{array}
\label{trivalt_rep}
\end{align}
%]

%%%%%%%%%%%%%%%%%%%%%%%%%%%%%%%%%%%%%%%%%%%%%%%%%%%%%%%%%%%%%%%%%%%%%%

\subsection{Symmetric trivial NCSR and the noncommutative sphere}
\label{ch_sym}

The situation is even simpler if $\rho$ is an even function with one
(quadratic) minima at $\rho(0)=0$.  In this case
%[
\begin{align}
\rho(-x)&=\rho(x) &
\tau(-x) &= -\tau(x) &
\omega(x) &= -\tfrac12 x
\end{align}
%]
Thus we have the much simpler
results:
%[
\begin{align}
N_k &= N_0 - \tfrac12\varepsilon k
&
X_k &= X_0 - \tfrac12\varepsilon k
&
[N_0,a_\pm] &= \mp\tfrac{\varepsilon}2 a_\pm
\notag
\\
[N_0,b_\pm] &= \mp\tfrac{\varepsilon}2 b_\pm
&
[X_0,a_\pm] &= \pm\tfrac{\varepsilon}2 a_\pm
&
[X_0,b_\pm] &= \mp\tfrac{\varepsilon}2 b_\pm
\end{align}
%]

\vspace{1 em}

Finally we give the connection between this and \cite{Gratus6}. If we
let $\rho(z)=z^2$ then $\tau(x)=-z$ and from 
(\ref{triv_com_NX}),(\ref{triv_prod_aa}) and (\ref{triv_com_ab}) we
deduce that $[a_\pm,b_\pm]=0$ and $[a_-,a_+]=[b_-,b_+]=\varepsilon$.
So $\Balg$ is simply the the product of two Heisenberg-Weyl algebras.

If we let $J_0:=X_0+\tfrac12\varepsilon$, $J_\pm=X_\pm$ and
$K_0:=-N_0$ we obtain the algebra given by the Jordan-Schwinger
representation of $su(2)$.  This is the starting point for the
analysis of vectors and spinors on the noncommutative sphere. To
obtain the representation in that article we relabel the vectors
$\rvec(n',m')$ where $n'=-\tfrac12(n-1)$ and $m'=m-\tfrac12(n-1)$.

%%%%%%%%%%%%%%%%%%%%%%%%%%%%%%%%%%%%%%%%%%%%%%%%%%%%%%%%%%%%%%%%%%%%%%
%%%%%%%%%%%%%%%%%%%%%%%%%%%%%%%%%%%%%%%%%%%%%%%%%%%%%%%%%%%%%%%%%%%%%%

\section{Vector module over NCSR}
\label{ch_vec}

In this section we still assume that $\rho$ is trivial that is it obeys
the constraints (\ref{rev_rho_asum}) and (\ref{triv_rho_asum}).
At the end of the last section we saw that the algebra $\Balg$ was,
for the special $\rho(z)=z^2$, simply the product of two
Heisenberg-Weil algebras. In \cite{Gratus6} we used this to produce
analogues of Vector and spinor fields for the noncommutative sphere.
In this section we repeat the process for NCSR.

\subsection{The algebra $(\Palg,\mu)$}
\label{ch_Psi}

Given $R,\varepsilon\in\Real$ and $C^1$ function
$\rho:\Real\mapsto\Real$ we define the non-associative algebra
$(\Palg,\mu)$.  Since this is a non-associative algebra we write the
product $\mu$ explicitly.  $\Palg$ is the set
%[
\begin{align}
\Palg &= \left\{
\xi=\sum_{r,s} a_\pm{}^r b_\pm{}^s \xi_{rs}(X_0)
\ \bigg|\
\xi_{rs}\in\Xdom
\text{ finite sum}  \right\}
\end{align}
%]
We define the non-associative product
$\mu:\Palg\times\Palg\mapsto\Palg$ as follows: Given two elements
$\xi,\zeta\in\Palg$, these may also be considered elements of
$\Balg$. We write the element $\xi\zeta\in\Balg$ in the form
(\ref{triv_B_gen_el}). Now make the identity $N_0\to R$ to produce an
element in $\Palg$ called $\mu(\xi,\zeta)$.  In other words
$(\Palg,\mu)=\Balg\big/{\cal I}(N_0-R)$ where we quotient $\Balg$ on
the right by the ideal generated by $N_0-R$.

We decompose $\Palg$ as:
%[
\begin{align}
\Palg = \bigoplus_{r\in\Intg} \Pmod{r} &
&
\text{where}&
&
\xi\in\Pmod{r}&\iff N_0 \xi = \xi N_{r}
\end{align}
%]
or equivalently
%[
\begin{align}
\xi\in\Pmod{r}&\iff 
\xi=\sum_m a_\pm^{r+m} b_\pm^{r-m} f^\xi_{m}(X_0)
\end{align}
%]
The set $\Pmod{0}$ is a subalgebra of $\Palg$, and it is equivalent to
$\Ralg$. All other sets $\Pmod{r}$ are
right modules over $\Ralg$. We wish to identify
these modules as vector fields, covectors fields, spinor fields, etc.
In this article we only interpret $\Pmod{2}$ and $\Pmod{-2}$ as the
space of vector and covector fields respectively. This is done in the
following section.

The relationship between our five algebras is given by:
%[
\begin{align}
\setlength{\unitlength}{8 mm}
\begin{picture}(7.5,2.5)
\put(0,2){\makebox(0,0){$\Calg$}}
\put(2,2){\makebox(0,0){$\NCalg$}}
\put(4,2){\makebox(0,0){$\Nalg$}}
\put(4,0){\makebox(0,0){$\Ralg$}}
\put(6,2){\makebox(0,0){$\Balg$}}
\put(6,0){\makebox(0,0){$\Palg$}}
\put(1,2){\makebox(0,0){$\hookrightarrow$}}
\put(5,2){\makebox(0,0){$\hookrightarrow$}}
\put(5,0){\makebox(0,0){$\hookrightarrow$}}
\put(2.5,2){\vector(1,0){1}}
\put(0,1.5){\vector(2,-1){3}}
\put(4,1.5){\vector(0,-1){1}}
\put(6,1.5){\vector(0,-1){1}}
\put(3,2.2){\makebox(0,0){${}^{q_1}$}}
\put(1.5,1){\makebox(0,0){${}^{q_2}$}}
\put(4.3,1){\makebox(0,0){${}^{q_3}$}}
\put(6.3,1){\makebox(0,0){${}^{q_4}$}}
\end{picture}
\label{Psi_com_diag}
\end{align}
%]
where the hooked arrows refer to the natural embedding $q_1$, $q_2$,
$q_3$, are given in lemma \ref{lm_rev_com_diag} and $q_4$ is 
the quotient $N_0-R\sim 0$ one the right.

%%%%%%%%%%%%%%%%%%%%%%%%%%%%%%%%%%%%%%%%%%%%%%%%%%%%%%%%%%%%%%%%%%%%%%
%%%%%%%%%%%%%%%%%%%%%%%%%%%%%%%%%%%%%%%%%%%%%%%%%%%%%%%%%%%%%%%%%%%%%%
%%%%%%%%%%%%%%%%%%%%%%%%%%%%%%%%%%%%%%%%%%%%%%%%%%%%%%%%%%%%%%%%%%%%%%
%%%%%%%%%%%%%%%%%%%%%%%%%%%%%%%%%%%%%%%%%%%%%%%%%%%%%%%%%%%%%%%%%%%%%%

\subsection{One forms over $\Ralg$}
\label{ch_onef}

We wish to define the space of one forms
$\Omega^1(\Ralg)$ and the exterior derivative
$d:\Ralg\mapsto\Omega^1(\Ralg)$. To do this we say that 
$\Omega^1(\Ralg)$ is a right module over $\Ralg$ which is spanned by
$\{\xi_0,\xi_+,\xi_-\}$ (which are not independent), such that
%[
\begin{align}
\xi_0 &= dX_0 &
\xi_+ &= dX_+ &
\xi_- &= dX_- 
\label{onef_xi}
\end{align}
%]
and there is a formula for $d$ which is consistent with (\ref{onef_xi})
and reduces to (\ref{poi_df}) in the limit $\varepsilon\to 0$.

Unfortunately, like most problems with quantisation, this procedure
is not unique. We are free to choose (\ref{poi_df2}) instead of
(\ref{poi_df}), we can choose the ordering of the elements, and we can
always add a random term which vanishes when $\varepsilon=0$.

We shall choose
\begin{align}
df 
&= 
\xi_+ X_-
\varepsilon^{-1} 
(\rho(R)-\rho(X_0))^{-1}
[X_0,f]
-
\xi_0 X_-
\varepsilon^{-1}(\rho(R)-\rho(X_0))^{-1}
[X_+,f]
\notag\\
&-
\tfrac12
\xi_0 X_-
\varepsilon^{-1}(\rho(R)-\rho(X_0+\varepsilon))^{-1}
\left(
[\rho(X_0+\varepsilon),f] - 
\varepsilon^{-1}(\rho(X_0+\varepsilon)-\rho(X_0))
\right)
\label{onef_df}
\end{align}
%]
where $\xi_+$ and $\xi_0$ are independent elements of $\Omega^1$.
We can see that $d$ is no longer a derivative for $\varepsilon\ne0$
but that
%[
\begin{align}
d(fh)=d(f)h+fd(h)+O(\varepsilon)
\end{align}
%]
Clearly $\xi_0=dX_0$ and $\xi_+=dX_+$. To be consistent with 
(\ref{onef_xi}) we let $\xi_-=dX_-$ giving
%[
\begin{align}
\xi_- 
&= 
-\xi_+
(\rho(R)-\rho(X_0+\varepsilon))^{-1}
X_-^2
-
\tfrac12
\xi_0
\varepsilon^{-1}
\frac{(\rho(X_0+2\varepsilon) -\rho(X_0))}
{(\rho(R)-\rho(X_0+\varepsilon))}
X_-
\label{onef_xim}
\\
\intertext{and rearranging}
\xi_+
&= 
-\xi_-
(\rho(R)-\rho(X_0))^{-1}
X_+^2
-
\tfrac12
\xi_0
\varepsilon^{-1}
\frac{(\rho(X_0+2\varepsilon) -\rho(X_0))}
{(\rho(R)-\rho(X_0+\varepsilon))}
X_+
\label{onef_xip}
\end{align}
%]
We can now write $df$ in terms of $\xi_0$ and $\xi_-$ giving
%[
\begin{align}
df &= 
-\xi_-
\varepsilon^{-1}(\rho(R)-\rho(X_0))^{-1}
X_+[X_0,f]
-
\xi_0
\varepsilon^{-1}(\rho(R)-\rho(X_0+\varepsilon))^{-1}
[X_-,f]X_+
\notag\\
&+
\tfrac12
\xi_0
\varepsilon^{-1}(\rho(R)-\rho(X_0+\varepsilon))^{-1}
\left(
[\rho(X_0+\varepsilon),f] - 
\varepsilon^{-1}(\rho(X_0+2\varepsilon)-\rho(X_0+\varepsilon))[X_0,f]
\right)
\label{onef_df2}
\end{align}
%]
We have chosen the final term in (\ref{onef_df}) in order to have the
maximum similarity between (\ref{onef_df}) and (\ref{onef_df2})

\vspace{1 em}

Both $\Omega^1(\Ralg)$ and $\Pmod{-2}$ are right modules
over $\Ralg$. We now make the identification
$\Omega^1(\Ralg)=\Pmod{-2}$ via the definitions
%[
\begin{align}
\xi_+ = a_-^2 \text{ and }
\xi_- = b_-^2
\end{align}
%]
and $\xi_0$ is given by (\ref{onef_xip}).

We can now construct the noncommutative analogue of the tangent bundle
$T\man$, a subset of which, we identify as $\Pmod{2}$. Given an
element $\xi\in\Pmod{2}$ we define the noncommutative analogue of the
vector field as the function $V_\xi:\Ralg\mapsto\Ralg$ given by
$V_\xi(f) = \mu(\xi,df)$.  We can see that for $\varepsilon\ne0$ then
$X$ does not obey Leibniz rule. As sated in the introduction this must
be the case since $T\man$ defined here is a right module.

By extending the results of \cite{Gratus6},
we can also interpret the modules $\Pmod{1}$ and $\Pmod{-1}$ as
spinors over $\Ralg$. This would mean the module
$\Pmod{r}\oplus\Pmod{-r}$ is the module of spin $r/2$ fields. 
This approach for spinors differs from the standard approach for
noncommutative geometry, using the supersymmetric group $SU(2|1)$
(for example \cite{Gross1}).

%%%%%%%%%%%%%%%%%%%%%%%%%%%%%%%%%%%%%%%%%%%%%%%%%%%%%%%%%%%%%%%%%%%%%%

\subsection{Derivatives of $\Ralg$ are inner}
\label{ch_inn}

As stated in the introduction, if we identify $T\man$ with the space
of derivations (obeying Leibniz) then it will not form a module over
$\Ralg$. This is because, as we shall show here, all derivations are
inner and it is easy to show that these do not form a module.

\begin{theorem}
\label{th_inner}
All Leibniz derivations on $\Ralg$ are inner.
\end{theorem}

\begin{proof}
Let $\xi:\Ralg\mapsto\Ralg$ be a derivations.  That is for any two
function $f,g\in\Ralg$, $\xi(fg)=\xi(f)g + f\xi(g)$. We are required to
show that there exists an $f\in\Ralg$ such that $\xi=\ad{f}$. Since
$\xi$ is a derivation it is only necessary to show that
$\xi(X_0)=\ad{f}(X_0)$ and $\xi(X_\pm)=\ad{f}(X_\pm)$. All other
functions can be derived from these.

Expanding $\xi(X_0)$ in terms of the eigenstates of $\ad{X_0}$,
written in normal form, we have
%[
\begin{align}
\xi(X_0) &= \sum_{r=-\infty}^\infty X_\pm^r p_r(X_0)
\notag
\end{align}
%]
For some set of functions $p_r$. Let
%[
\begin{align}
\wh f= \sum_{r=-\infty,r\ne 0}^\infty \tfrac1r X_\pm^r p_r(X_0)
\notag
\end{align}
%]
Then
%[
\begin{align}
\xi(X_0) &= [\wh f,X_0] + p_0(X_0)
\notag
\end{align}
%]
Since $\xi$ and $\ad{\wh f}$ are both derivative then this formula
extends to any function of $X_0$ as
%[
\begin{align}
\xi(h(X_0)) &= [\wh f,h(X_0)] + h'(X_0)p_0(X_0)
\notag
\end{align}
%]

Let us define
%[
\begin{align}
g_\pm &= \xi(X_\pm) - [\wh f,X_\pm]
\notag
\end{align}
%]
Taking the $\xi$ derivative of (\ref{rev_com_X0Xpm}) gives
%[
\begin{align}
[\xi(X_0),X_\pm] + [X_0,X_\pm] &= \pm\varepsilon\xi(X_\pm)
\notag
\end{align}
%]
Expanding and substituting the above expressions give
%[
\begin{align}
[X_0,g_\pm] &= \pm\varepsilon g_\pm + [X_\pm,p_0(X_0)]
\notag
\end{align}
%]
By expanding in normal form $g_+$ we see that
%[
\begin{align}
[p_0(X_0),X_+] &= \varepsilon \sum_{r=-\infty}^\infty (1-r) X_\pm^r
g_r(X_0) 
\notag
\end{align}
%]
Which is only consistent if $p_0=0$. This means that
$[X_0,g_+] = \varepsilon g_+$. This means we can write
$g_+=X_+(g(X_0+\varepsilon)-g(X_0))$ for some function $g$.
Let $f=\wh f + g(X_0)$. Then 
%[
\begin{align}
\xi(X_0) &= [f,X_0] & 
\xi(X_+) &= [f,X_+] &
\notag
\end{align}
%]
Now taking the derivative of (\ref{rev_Ra_XX}) and expanding gives
%[
\begin{align}
\xi(X_+)X_- + X_+\xi(X_-) &= \xi(\rho(R)-\rho(X_0)) 
\notag
\\
[\wh f,X_+]X_- + g_+X_- + X_+[\wh f,X_-] + X_+g_- &=
-[\wh f,\rho(X_0)] 
\notag
\\
g_+ X_- + X_+ g_- &= 0
\notag
\\
X_+(g(X_0+\varepsilon)-g(X_0))X_- + X_+ g_- &= 0 
\notag
\end{align}
%]
which gives $g_-=[g(X_0),X_-]$. Which implies $\xi(X_-)=[f,X_-]$.
\end{proof}

%%%%%%%%%%%%%%%%%%%%%%%%%%%%%%%%%%%%%%%%%%%%%%%%%%%%%%%%%%%%%%%%%%%%%%
%%%%%%%%%%%%%%%%%%%%%%%%%%%%%%%%%%%%%%%%%%%%%%%%%%%%%%%%%%%%%%%%%%%%%%
%%%%%%%%%%%%%%%%%%%%%%%%%%%%%%%%%%%%%%%%%%%%%%%%%%%%%%%%%%%%%%%%%%%%%%

\section{Representations of non trivial NCSR}
\label{ch_notriv}

Having established how to represent trivial NCSR we know turn our
attention to how to represent non-trivial NCSR. Thus we assume that
$\rho$ obeys the conditions given by (\ref{rev_rho_asum}), and that it
has at least two local minima.

There are problems associated with the construction of the
multi-topology representations of $\rho$. Within each topology
interval the representation of $\rho$ is similar to the representation
of a trivial $\rho$. It is at the boundaries of the topology intervals
that the problems occur. This is not surprising since it is here that
the topology changes occurs.

In this section we do not try to construct the algebra $\Balg$ but
simply try to construct a Hilbert space $\Gspace$ and the hopping
operators. An example $\Gspace$ together with over hopping operators
is given in figure \ref{fig_notriv}.

In section \ref{ch_G0} we construct the Hilbert space $\Gspace_0$, and
in section \ref{ch_Hop} we construct the hopping operators. As
indicated in section \ref{ch_comp}, these do not form a true
representation of $\Nalg$ and it is necessary to compromise.  There
are at least four possible compromises. Two of them use over-hopping
operators and two use under-hopping operators. Unlike the case of a
trivial $\rho$, in general there is no isomorphism between the
over-hopping operators and the under-hopping operators.

%%%%%%%%%%%%%%%%%%%%%%%%%%%%%%%%%%%%%%%%%%%%%%%%%%%%%%%%%%%%%%%%%%%%%%

\subsection{The Hilbert space $\Gspace_0$}
\label{ch_G0}

We construct first the Hilbert space $\Gspace_0$. For two of the four
choices for representations we use $\Gspace_0$ directly for the others
we have to modify $\Gspace_0$.  Given a $\rho$ which obeys
(\ref{rev_rho_asum}) and an $\varepsilon>0$, the Hilbert space
$\Gspace_0=\Gspace_0(\rho,\varepsilon)$ is defined as follows:

Let $\{J_s\subset\Real,\ s\in\natnum\}$ be a set of intervals obeying
(\ref{rep_def_J}). The end points of the interval $J_s$ are given by
points $z\smin{s}$ and $z\smax{s}$ so $J_s=\{z\in\Real\ |\
z\smin{s}\le z\le z\smax{s}\}$.  Define the function
$\pi:\natnum\mapsto\natnum$ as: Given $J_s$ then $J_{\pi(s)}$ is the
smallest $J_t\ne J_s$ such that $J_s\subset J_t$.

Recall that since $\T$ is a finite set there exists an element
$t_\infty$ such that $\piT(t_\infty)=\infty$. This is the unbounded
topology interval.  
For each $s\in\S$ we define integers $n_s$, $m\smin{s}$ and
$m\smax{s}$ as:
Choose an $s_0$ such that
$J_{s_0}\in\topI{t_\infty}$.  
Let $n_{s_0}=0$ and $m\smin{s_0}=0$.
Now for all $J_s\in\topI{t_\infty}$ let $n_{\pi(s)}=n_s+1$ and
$m\smin{s}=0$. For the other $s$ let $n_s=n_{\pi(s)}-1$ and
%[
\begin{align}
m\smin{s}=\left\lfloor 
\frac{z\smin{s}-z\smin{\pi(s)}}{\varepsilon} \right\rfloor
+m\smin{\pi(s)}
\qquad\qquad
m\smax{s}=m\smin{s}+\frac{|J|}{\varepsilon} -1
\end{align}
%]
For each $s\in\S$ let $V_s$ be the vector space spanned by the basis
vectors $\{\rvec(n_s,m)\}$ where $m\in\Intg$, $m\smin{s}\le m\le
m\smax{s}$.

Finally we define $\Gspace_0$ as $\Gspace_0=\oplus_{s\in\natnum}
V_s$. An example of a $\Gspace_0$ is given in figure \ref{fig_notriv}.

\vspace{1 em}

\begin{figure}[t]
\begin{center}
\epsfysize 5 cm
\centerline{\epsffile{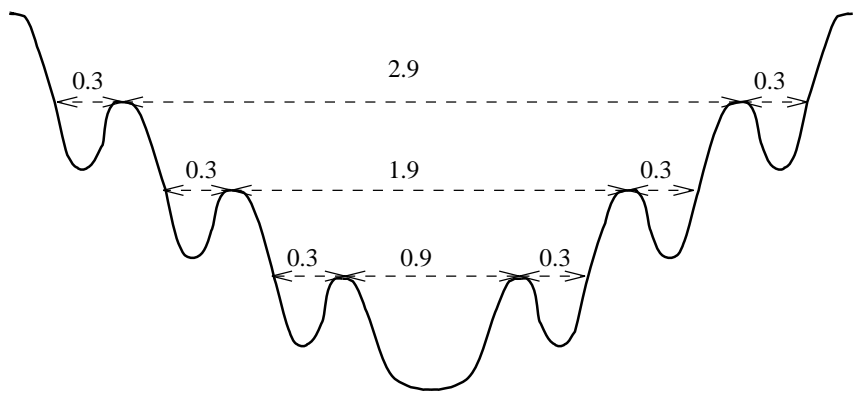}}
\caption{Example of a $\rho$ for which no \Gspace-representation exists
($\varepsilon=1$)}
\label{fig_no_Grat}
\end{center}
\vspace{10em}
\end{figure}

We have required that $\rho$ has only a finite number of maxima since
otherwise it is possible to construct $\rho$ such that no, or only a
finite number of $J_s$ exist. See figure \ref{fig_no_Grat}.

The choice of $s$ for each $J_s$ is arbitrary in the construction of
$\Gspace_0$ as is the initial $s_0$, but once these are set, that
fixes $n$ and $m$ for each basis vector $\rvec(n,m)$.

%%%%%%%%%%%%%%%%%%%%%%%%%%%%%%%%%%%%%%%%%%%%%%%%%%%%%%%%%%%%%%%%%%%%%%

\subsection{Hopping operators}
\label{ch_Hop}

For each $\rvec(n,m)\in\Gspace_0$ there is a unique $s\in\natnum$ such
that $\rvec(n,m)\in V_s$. The effect of the diagonal operators are
given by
%[
\begin{align}
f(X_0,N_0)=
f(z\smin{s}+\varepsilon m-\varepsilon m\smin{s},z\smin{s})
\rvec(n,m) 
\label{Hop_XN}
\end{align}
%]
For each $s\in\natnum$ and $m\smin{s}\le m\le m\smax{s}$ let
%[
\begin{align}
D_s &= \min\{\rho(z),z\in J_s\}
\label{Hop_def_Ds}
\\
(\SS(s,m))^2 &= \rho(x\smin{s}-\varepsilon m\smin{s}+\varepsilon m)
- D_s
\label{Hop_def_ss2}
\\
C_s &= \SS(s,m\smin{s} ) = -
\SS(s,m\smax{s}  +1) = 
\left(\rho(x\smin{s}) - D_s \right)^\scrhalf
> 0
\label{Hop_def_Cs}
\end{align}
%]
The sign of $\SS(s,m)$ is arbitrary except for $m=m\smin{s}$ and
$m=m\smax{s}+1$.  The construction of $\Gspace_0$ guarantees that if
$\rvec(n,m)\in V_s$ then $\rvec(n+1,m)\in\Gspace_0$ and
$\rvec(n+1,m+1)\in\Gspace_0$.  In contrast to the trivial case, the
over-hopping operators and under-hopping operators are no longer
equivalent.  The over-hopping operators are given by
%[
\begin{eqnarray}
\begin{aligned}
a_+\rvec(n,m) &= (C_s-\SS(s,m+1))^\scrhalf\rvec(n+1,m+1) 
\\
a_-\rvec(n,m) &= 
\begin{cases}
(C_t-\SS(t,m))^\scrhalf\rvec(n-1,m-1) 
& \text{ if } \rvec(n-1,m-1)\in V_t 
\text{ for some }V_t\subset\Gspace_0 \\
\makebox[9 em]{0} &\text{ if } \rvec(n-1,m-1)\not\in\Gspace_0
\end{cases}
\\
b_+\rvec(n,m) &= (C_s+\SS(s,m))^\scrhalf\rvec(n+1,m) 
\\
b_-\rvec(n,m) &= 
\begin{cases}
(C_t+\SS(t,m))^\scrhalf\rvec(n-1,m) 
& \text{ if } \rvec(n-1,m)\in V_t 
\text{ for some }V_t\subset\Gspace_0 \\
\makebox[9 em]{0} &\text{ if } \rvec(n-1,m)\not\in\Gspace_0
\end{cases}
\end{aligned}
\label{Hop_over_abpm}
\end{eqnarray}
%]
The under-hopping operators are given (with $t=\pi(s)$) by
%[
\begin{eqnarray}
\begin{aligned}
a'_+\rvec(n,m) &= (C_t-\SS(t,m+1))^\scrhalf\rvec(n+1,m+1) 
\\
a'_-\rvec(n,m) &= 
\begin{cases} 
(C_s-\SS(s,m))^\scrhalf\rvec(n-1,m-1) 
& \text{ if } \rvec(n-1,m-1)\in\Gspace_0
\\
\makebox[9 em]{0} &\text{ if } \rvec(n-1,m-1)\not\in\Gspace_0 
\end{cases}
\\
b'_+\rvec(n,m) &= (C_t+\SS(t,m))^\scrhalf\rvec(n+1,m) 
\\
b'_-\rvec(n,m) &= 
\begin{cases} 
(C_s+\SS(s,m))^\scrhalf\rvec(n-1,m) 
& \text{ if } \rvec(n-1,m-1)\in\Gspace_0
\\
\makebox[9 em]{0} &\text{ if } \rvec(n-1,m)\not\in\Gspace_0 
\\
\end{cases}
\end{aligned}
\label{HO_under_abpm}
\end{eqnarray}
%]
In order to compare the hopping operators with the ladder operators,
we define the action of the ladder operators on $\Gspace_0$ as:
%[
\begin{equation}
\begin{aligned}
{}&
X_+\rvec(n,m) = 
(\rho(z\smin{s})-
\rho(z\smin{s}-\varepsilon m\smin{s}+\varepsilon m+\varepsilon)
)^\scrhalf\rvec(n,m+1)
\\
{}&
X_-\rvec(n,m) = 
(\rho(z\smin{s})-
\rho(z\smin{s}-\varepsilon m\smin{s}+\varepsilon m)
)^\scrhalf\rvec(n,m-1)
\end{aligned}
\label{Hop_Xpm}
\end{equation}
%]
These are consistent with the $(J_s,V_s)$ representation of $\Nalg$
given by (\ref{rep_stand_rep}).

\begin{lemma}
\label{lm_nontriv_Xpm}
The operators $N_0$ and $X_0$ are self-adjoint and the operators
$a_+,b_+,a'_+,b'_+$ are the adjoints of $a_-,b_-,a'_-,b'_-$
respectively. The points given by (\ref{trivrep_coords}) are well
placed since (\ref{Hop_XN}) satisfies (\ref{trivrep_inc_rhoN}) and
(\ref{trivrep_posn_X0}). The hopping operators are related to the
ladder operators in the following circumstances:
%[
\begin{align*}
&\text{if }
\rvec(n,m)\in V_s
\,,\ 
\rvec(n,m+1)\in V_s
\text{ then }
b_-a_+\rvec(n,m)
=
(C_s^2 - (\SS(s,m+1))^2)^\scrhalf\rvec(n,m+1)
=
X_+\rvec(n,m)
\\
&\text{if }
\rvec(n,m)\in V_s
\,,\ 
\rvec(n,m-1)\in V_s
\text{ then }
a_-b_+\rvec(n,m)
=
(C_s^2 - (\SS(s,m+1))^2)^\scrhalf\rvec(n,m-1)
=
X_-\rvec(n,m)
\\
&
\text{if }
\rvec(n,m)\in V_s
\,,\ 
\rvec(n,m+1)\not\in \Gspace_0
\text{ then }
b_-a_+\rvec(n,m)=0
\\
&
\text{if }
\rvec(n,m)\in V_s
\,,\ 
\rvec(n,m-1)\not\in \Gspace_0
\text{ then }
a_-b_+\rvec(n,m)=0
\\[0.5 em]
&\text{if }
\rvec(n,m)\in V_s
\,,\ 
\rvec(n,m+1)\in V_s
\,,\ 
\rvec(n-1,m)\in \Gspace_0
\text{ then }
a'_+b'_-\rvec(n,m)
=
(C_s^2 - (\SS(s,m+1))^2)^\scrhalf\rvec(n,m+1)
=
X_+\rvec(n,m)
\\
&\text{if }
\rvec(n,m)\in V_s
\,,\ 
\rvec(n,m-1)\in V_s
\,,\ 
\rvec(n-1,m-1)\in \Gspace_0
\text{ then }
b'_+a'_-\rvec(n,m)
=
(C_s^2 - (\SS(s,m))^2)^\scrhalf\rvec(n,m-1)
=
X_-\rvec(n,m)
\\
&\text{if }
\rvec(n,m)\in V_s
\,,\ 
\rvec(n,m+1)\not\in \Gspace_0
\text{ or }
\rvec(n,m)\in V_s
\,,\ 
\rvec(n-1,m)\not\in\Gspace_0 
\text{ then }
a'_+b'_-\rvec(n,m)
= 0
\\
&\text{if }
\rvec(n,m)\in V_s
\,,\ 
\rvec(n,m+1)\not\in \Gspace_0
\text{ or }
\rvec(n,m)\in V_s
\,,\ 
\rvec(n-1,m-1)\not\in\Gspace_0 
\text{ then }
b'_+a'_-\rvec(n,m)
= 0
\end{align*}
%]
Thus $(J_s,V_s)$ is a representation of $\Nalg$ using over-hopping
operators where $X_+=b_-a_+$ and $X_-=a_-b_+$ if
%[
\begin{align}
\dim V_s \ge \sum_{t\in\pi^{-1}\{s\}} (\dim V_t + 1)
\label{Hop_dim_over}
\end{align}
%]
On the other hand
$(J_s,V_s)$ is a representation of $\Nalg$ using under-hopping
operators where $X_+=a'_+b'_-$ and $X_-=b'_+a'_-$ if
%[
\begin{align}
\dim V_s = \left(\sum_{t\in\pi^{-1}\{s\}} \dim V_t\right) + 1
\label{Hop_dim_under}
\end{align}
%]
\end{lemma}

\begin{proof}
The first and second parts follow from direct substitution.
The conditions for $\Gspace_0$ to be a representation of $\Nalg$ 
follows from the requirement that $X_+\rvec(n,m)=0$ if and only if
$n=n_s$ and $m=m\smax{s}$ for some $s\in\natnum$. An example which
obeys (\ref{Hop_dim_over}) is given in figure \ref{fig_notriv}.
\end{proof}

\begin{figure}[t]
\begin{center}
\epsfysize 5 cm
\centerline{\epsffile{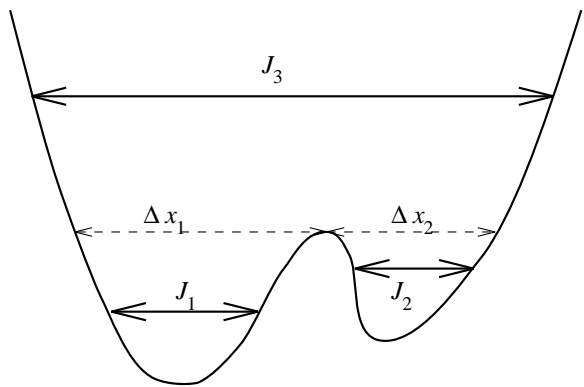}}
\caption{An example of a $\rho$ containing a maxima}
\label{fig_ignor}
\end{center}
\vspace{10em}
\end{figure}

Clearly for a given $\rho$ and $\varepsilon$ both (\ref{Hop_dim_over})
and (\ref{Hop_dim_under}) cannot be satisfied. To see which, if either, is
satisfied consider figure
(\ref{fig_ignor}).  Here $J_1$, $J_2$, $J_3$ are intervals satisfying
(\ref{rep_def_J}). Thus
%[
\begin{align}
|J_1|&=\varepsilon \lfloor \Delta x_1/\varepsilon \rfloor, &
|J_2|&=\varepsilon \lfloor \Delta x_2/\varepsilon \rfloor, &
|J_3|&=\varepsilon \lceil (\Delta x_1+\Delta x_2)/\varepsilon \rceil
\end{align}
%]
where $\lfloor x\rfloor$ and $\lceil x\rceil$ are the nearest integer
below and above $x$ respectively.  As a result either
$|J_3|=|J_1|+|J_2|+\varepsilon$ obeying (\ref{Hop_dim_under}) or
$|J_3|=|J_1|+|J_2|+2\varepsilon$ obeying (\ref{Hop_dim_over}).  If
$\rho$ has just one maxima then we can choose over-hopping or
under-hopping operators appropriately. However if $\rho$ has more than
one maxima it may be impossible to satisfy either (\ref{Hop_dim_under})
for all $s\in\natnum$ or (\ref{Hop_dim_over}) for all $s\in\natnum$.

%%%%%%%%%%%%%%%%%%%%%%%%%%%%%%%%%%%%%%%%%%%%%%%%%%%%%%%%%%%%%%%%%%%%%%

\subsection{Four compromises}
\label{ch_comp}

For each $s\in\natnum$ we have a classical surface $\man_s$ given by
%[
\begin{align}
\man_s = \{(x^1,x^2,x^3)\in\Real^3\ |\
(x^1)^2+(x^2)^2=\rho(R)-\rho(x^3), x^3\in J_s \}
\end{align}
%]
We would like to have $\Gspace_0$ a unitary representation of $\Nalg$
and each $V_s$ as the standard unitary representation of
$\Ralg(\rho,\varepsilon,R)$, where $R=z\smin{s}$ and thus the
noncommutative analogue of $\man_s$.  As noted, in general, either
(\ref{Hop_dim_under}) or (\ref{Hop_dim_over}) is not satisfied.  It is
therefore necessary either to accept that there exists some intervals
$J_s$ where we do not have a true standard unitary of $\Nalg$, or to
alter $\Gspace_0$ in some way. There are at least four possible
compromises, two of which use over-hopping operators and two of which
use under-hopping operators. The disadvantages of each compromise are
summarised in the following table:

\noindent
\begin{tabular}{|l|c|l|l|l|l|}
\hline
Comprise &
Hop op &
Uses $\Gspace_0$ &
Rep of $\Nalg$ &
Obeys (\ref{trivrep_inc_rhoN}) &
Obeys (\ref{trivrep_posn_X0})
\\
\hline
1. Merging $V_s$ &
over &
yes &
no &
yes &
yes 
\\
2. Splitting $V_s$ &
under &
yes &
no &
yes &
yes 
\\
3. Removing $V_s$ &
over &
no &
yes &
yes &
no 
\\
4. Add a vector &
under &
no &
yes &
yes &
yes 
\\
\hline
\end{tabular}

For the following we assume that $\rho$ has only one maxima and two
minima. However we can see how this generalises for any finite number
of maxima.

%%%%%%%%%%%%%%%%%%%%%%%%%%%%%%%%%%%%%%%%%%%%%%%%%%

\subsubsection{Compromise 1: Merging of $V_s$s}

\begin{figure}[t]
\begin{center}
\setlength{\unitlength}{0.5 mm}
\begin{picture}(300,100)
\put(0,0){\epsfysize 100 \unitlength
{\epsffile{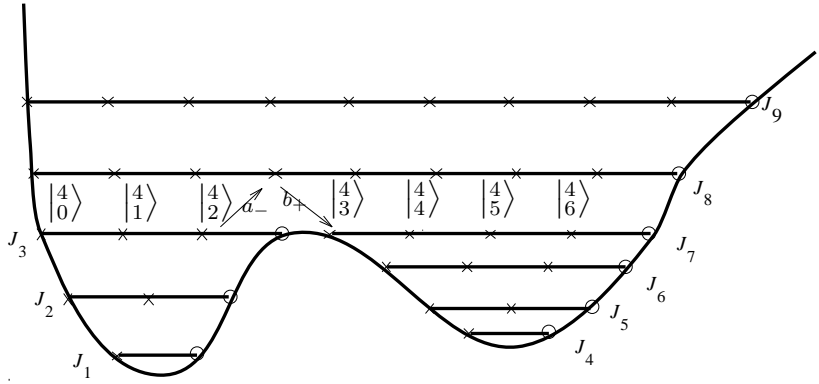}}}
%\graphpaper(0,0)(300,100)
\put(0,0){\setcounter{nn}{0}}
\multiput(10,42)(20,0){3}
{\makebox(0,0)[bl]{$\rvec(4,{\arabic{nn}})$\addtocounter{nn}{1}}}
\multiput(85,43)(20,0){4}
{\makebox(0,0)[bl]{$\rvec(4,{\arabic{nn}})$\addtocounter{nn}{1}}}
\put(62,47){\makebox(0,0)[tl]{${}^{a_-}$}}
\put(80,50){\makebox(0,0)[tr]{${}^{b_+}$}}
\end{picture}
\end{center}
\caption{Merging of $V_3$ and $V_7$}
\label{fig_merge}
\vspace{10em}
\end{figure}

We use the Hilbert space $\Gspace=\Gspace_0$ and the over-hopping
operators.  Consider figure \ref{fig_merge} where $\dim V_8=\dim V_3 +
\dim V_7 +1$ and hence violates (\ref{Hop_dim_over}).  If we let
$V_3=\spanrm\{\rvec(4,0),\rvec(4,1),\rvec(4,2)\}$ and
$V_7=\spanrm\{\rvec(4,3),\rvec(4,4),\rvec(4,5),\rvec(4,6)\}$ then
%[
\begin{align}
X_+ \rvec(4,2) &= (4C_3C_7)^\scrhalf\rvec(4,3) &
X_- \rvec(4,3) &= (4C_3C_7)^\scrhalf\rvec(4,2) 
\notag
\end{align}
%]
If $V_3$ and $V_7$ were both representations of $\Nalg$ then
$X_+\rvec(4,2)=X_-\rvec(4,3)=0$. We say that $V_3$ and $V_7$ have
merged. 
We also see that 
%[
\begin{align}
\lvec(4,3)[N_0,X_+]\rvec(4,2) = (4C_3C_7)^\scrhalf (z\smin{7}-z\smin{3}) \ne 0
\end{align}
%]
thus violating (\ref{rev_com_N0Xpm}).  For all other $V_s$ except
$V_3$ and $V_7$ there is a valid representation of $\Nalg$.

An interpretation of this is that the two surfaces $\man_3$ and
$\man_7$ are less apart than the distance of the parameter
$\varepsilon$. Since we interpret the noncommutative versions of
$\man_3$ and $\man_7$ as somehow smearing out the classical surfaces
then this causes $\man_3$ and $\man_7$ to interact.

%%%%%%%%%%%%%%%%%%%%%%%%%%%%%%%%%%%%%%%%%%%%%%%%%%

\subsubsection{Compromise 2: Splitting $V_s$s}

\begin{figure}[t]
\begin{center}
\setlength{\unitlength}{0.5 mm}
\begin{picture}(300,100)
\put(0,0){\epsfysize 100 \unitlength
{\epsffile{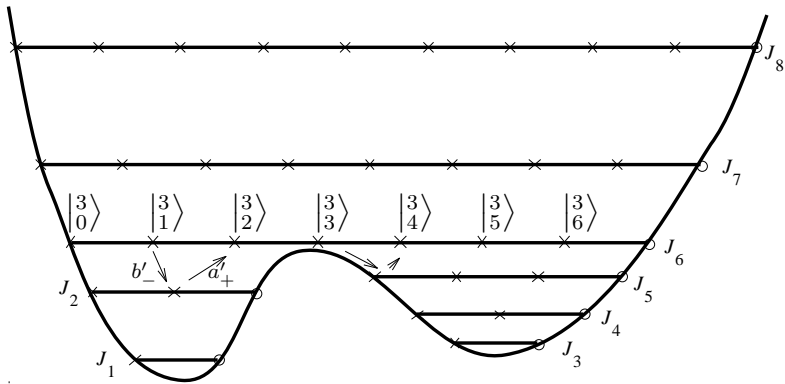}}}
%\graphpaper(0,0)(300,100)
\put(0,0){\setcounter{nn}{0}}
\multiput(15,40)(22,0){7}{\makebox(0,0)[bl]{$\rvec(3,{\arabic{nn}})$\addtocounter{nn}{1}}}
\put(40,32){\makebox(0,0)[tr]{${}^{b'_-}$}}
\put(53,32){\makebox(0,0)[tl]{${}^{a'_+}$}}
\end{picture}
\caption{Splitting of $V_8$}
\label{fig_split}
\end{center}
\vspace{10em}
\end{figure}

We use the Hilbert space $\Gspace=\Gspace_0$ and the under-hopping
operators. Consider figure \ref{fig_split} where $\dim V_6=\dim V_2 +
\dim V_5 +2$ and hence violates (\ref{Hop_dim_under}).  If we let
$V_6=\spanrm\{\rvec(3,0),\rvec(3,1),\ldots,\rvec(3,6)\}$ then
%[
\begin{align}
X_+\rvec(3,2)=X_-\rvec(3,3)=0
\end{align}
%]
thus violating (\ref{rev_Na_f}).

This has the opposite interpretation to the compromise above. 
Because the classical surface $\man_6$ has a small waist its
noncommutative analogue is equivalent to two surfaces.

%%%%%%%%%%%%%%%%%%%%%%%%%%%%%%%%%%%%%%%%%%%%%%%%%%

\subsubsection{Compromise 3: Removing $V_s$}

Here we use the over-hopping operators, but modify $\Gspace_0$ if
there exist a $V_s$ which violates (\ref{Hop_dim_over}).

Let $\Gspace=\oplus_{s\in\S}$ where $\S\in\natnum$ is the set of all
$V_s$ which obey (\ref{Hop_dim_over}), i.e. we remove from
$\Gspace_0$ all $V_s$ which violate (\ref{Hop_dim_over}).
We then have to redefine the hopping operators. 

\begin{figure}[t]
\begin{center}
\setlength{\unitlength}{0.5 mm}
\begin{picture}(300,100)
\put(0,0){\epsfysize 100 \unitlength
{\epsffile{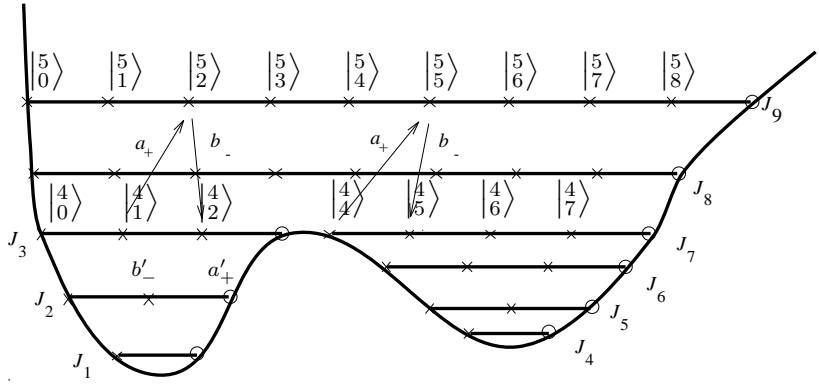}}}
%\graphpaper[10](0,0)(200,100)
\put(0,0){\setcounter{nn}{0}}
\multiput(10,42)(20,0){3}
{\makebox(0,0)[bl]{$\rvec(4,{\arabic{nn}})$\addtocounter{nn}{1}}}
\put(0,0){\setcounter{nn}{4}}
\multiput(85,43)(20,0){4}
{\makebox(0,0)[bl]{$\rvec(4,{\arabic{nn}})$\addtocounter{nn}{1}}}
\put(0,0){\setcounter{nn}{0}}
\multiput(5,77)(21,0){9}
{\makebox(0,0)[bl]{$\rvec(5,{\arabic{nn}})$\addtocounter{nn}{1}}}
\put(40,32){\makebox(0,0)[tr]{${}^{b'_-}$}}
\put(53,32){\makebox(0,0)[tl]{${}^{a'_+}$}}
\end{picture}
\end{center}
\caption{Removal of $V_7$}
\label{fig_rmline}
\vspace{10em}
\end{figure}

Consider figure \ref{fig_rmline} where $\dim V_8=\dim V_3 + \dim V_7
+1$ and hence violates (\ref{Hop_dim_over}).  We do not include $V_8$
from $\Gspace$. We have to relabel $\rvec(n,m)$ for all $V_s$ where
$J_s\subset J_8$.  If we let
$V_9=\spanrm\{\rvec(5,0),\rvec(5,1),\ldots,\rvec(5,8)\}$ then we must
make $V_3=\spanrm\{\rvec(4,0),\rvec(4,1),\rvec(4,2)\}$ and
$V_7=\spanrm\{\rvec(4,4),\rvec(4,5),\rvec(4,6),\rvec(4,7)\}$ so that
there is a vector missing between $V_3$ and $V_7$. Otherwise these
spaces would merge.  This requirement may mean that
(\ref{trivrep_posn_X0}) is violated as in figure \ref{fig_rmline}.

It is difficult to see how to interpret this compromise but at least
$\Gspace$ is a unitary representation of $\Nalg$ and all the $V_s$
that remain in $\Gspace$ are unitary representation of $\Ralg$.

%%%%%%%%%%%%%%%%%%%%%%%%%%%%%%%%%%%%%%%%%%%%%%%%%%

\subsubsection{Compromise 4: Adding an extra vector}

If we use the under-hopping operators so $X_+=a'_+b'_-$ and
$X_-=b'_+a'_-$ and modify $\Gspace_0$ if there exist a $V_s$ that
violates (\ref{Hop_dim_under}).

\begin{figure}[t]
\begin{center}
\setlength{\unitlength}{0.5 mm}
\begin{picture}(300,100)
\put(0,0){\epsfysize 100 \unitlength
{\epsffile{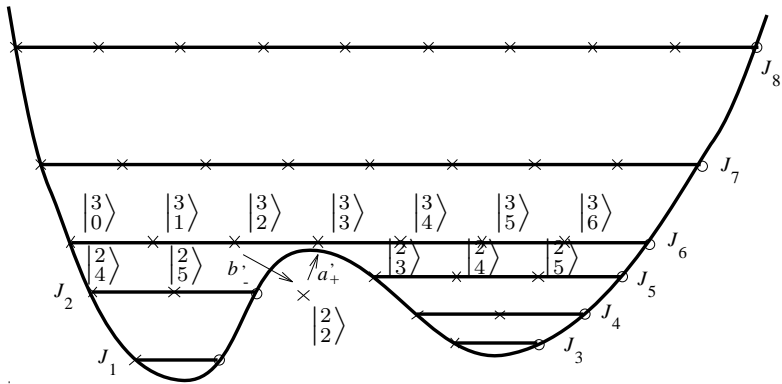}}}
%\graphpaper[10](0,0)(200,100)
\put(0,0){\setcounter{nn}{0}}
\multiput(19,40)(22,0){7}
{\makebox(0,0)[bl]{$\rvec(3,{\arabic{nn}})$\addtocounter{nn}{1}}}
\put(0,0){\setcounter{nn}{4}}
\multiput(20,26)(22,0){2}
{\makebox(0,0)[bl]{$\rvec(2,{\arabic{nn}})$\addtocounter{nn}{1}}}
\put(0,0){\setcounter{nn}{3}}
\multiput(100,28)(21,0){3}
{\makebox(0,0)[bl]{$\rvec(2,{\arabic{nn}})$\addtocounter{nn}{1}}}
\put(80,20){\makebox(0,0)[tl]{$\rvec(2,2)$}}
\end{picture}
\end{center}
\caption{Adding a point $V_6$}
\label{fig_addpt}
\vspace{10em}
\end{figure}

Consider figure \ref{fig_addpt} where $\dim V_6=\dim V_2 + \dim V_5
+2$ and hence violates (\ref{Hop_dim_under}).  Let
$V_6=\spanrm\{\rvec(3,0),\rvec(3,1),\ldots,\rvec(3,6)\}$,
$V_2=\spanrm\{\rvec(2,0),\rvec(2,1)\}$ and
$V_5=\spanrm\{\rvec(2,3),\rvec(2,4),\rvec(2,5)\}$. Now let
$\Gspace=\Gspace_0\oplus\spanrm\{\rvec(2,2)\}$. By enlarging $\Gspace$
we have avoided the problem in compromise 2.  We define $a_+$ and
$b_+$ on $\rvec(2,2)$ by considering $a_-\rvec(3,3)$ and
$b_-\rvec(3,2)$.

Again $\Gspace$ is a unitary representation of $\Nalg$ and all the $V_s$
that remain in $\Gspace$ are unitary representation of $\Ralg$.

%%%%%%%%%%%%%%%%%%%%%%%%%%%%%%%%%%%%%%%%%%%%%%%%%%%%%%%%%%%%%%%%%%%%%%

\subsection{Representation of the topology change}
\label{ch_rep_topch}

We say that a lattice representation $\Gspace(\rho,\varepsilon)$
reflects the topology of $\rho$ if 
\begin{jglist}
\item for all $t\in\T$ there exists $s\in\natnum$ such that
$J_s\in\topI{t}$.

\item for all $s\in\natnum$ there exists $t\in\T$ such that 
$J_s\in\topI{t}$. Either $J_{\pi(s)}\in\topI{t}$ or
$J_{\pi(s)}\in\topI{\piT(t)}$.

\item for all $t\in\T$ such that $\piT(t)\ne\infty$ then there exists
a unique $s\in\natnum$ such that $J_s\in\topI{t}$ and
$J_{\pi(s)}\in\topI{\piT(t)}$.
\end{jglist}
This definition gives the relationship between the maps
$\pi:\natnum\mapsto\natnum$  and $\piT:\T\mapsto\T\union\{\infty\}$. 

\begin{theorem}
For all four compromise there exist an $\varepsilon_0>0$ such that for
al $0<\varepsilon<\varepsilon_0$ there exists a lattice representation
of the NCSR which respects the topology of $\rho$.
\end{theorem}

\begin{proof}
Clearly for all $\varepsilon$ we can construct the space $\Gspace_0$
and hence the space $\Gspace$. If we let 
%[
\begin{align*}
\varepsilon_0 = \min_{t\in\T}
\bigg(
\sup_{I\in\topI{t}}(|I|)
-
\inf_{I\in\topI{t}}(|I|)
\bigg)
\end{align*}
%]
This guarantees that $\Gspace_0$ respects the topology of $\rho$. Thus
the result for  compromises 1, 2, and 4. For compromise 4 we note that
the added vectors do not correspond to any interval $J_s$ and thus are
not in a topology interval $\topI{t}$.

For compromise 3 it is necessary to replace $\varepsilon_0$ with
$\varepsilon_0/|\T|$ where $|\T|$ are the number of elements in
$\T$. This guarantees that even if we remove the maximum number of
spaces $V_s$ there is still a $J_s$ in each $\topI{t}$. 
\end{proof}

%%%%%%%%%%%%%%%%%%%%%%%%%%%%%%%%%%%%%%%%%%%%%%%%%%%%%%%%%%%%%%%%%%%%%%

\subsection{Operators for topology change}
\label{ch_topch}

For each $V_s$ let $L(V_s)$ be the set of operators on $V_s$
and let $\pi^{-1}(V_s)=\oplus_{t\in\pi^{-1}(s)}V_t$.
Let us define the linear operators
%[
\begin{equation}
\begin{aligned}
{}&\Aop_+,\Bop_+:L(\pi^{-1}(V_s))\mapsto L(V_s)
&{\qquad\qquad}&\Aop_-,\Bop_-:L(V_s)\mapsto L(\pi^{-1}(V_s))
\\
{}&\Aop_+(f) = a_+(a_-a_+)^{-1}f a_-
&{}&\Aop_-(f) = a_-f a_+(a_-a_+)^{-1}
\\
{}&\Bop_+(f) = b_+(b_-b_+)^{-1}f b_-
&{}&\Bop_-(f) = b_-f b_+(b_-b_+)^{-1}
\end{aligned}
\label{topch_def_ABpm}
\end{equation}
%]
\begin{lemma}
The maps ${\Aop_\pm},{\Bop_\pm}$ are well defined. $\Aop_+$ and
$\Bop_+$ are homomorphism i.e.  ${\Aop_+}(fg)={\Aop_+}(f){\Aop_+}(g)$.
Given $f\in L(\pi^{-1}(V_s)$ then
%[
\begin{align}
\Aop_-(\Aop_+(f))=\Bop_-(\Bop_+(f))=f
\label{topch_hom}
\end{align}
%]
We can write the elements in ${L(V_s)}$ as matrices using the
basis $\rvec(n,m)$.  As a matrix, ${\Aop_+}(f)$ has zeros in
the first row and column, and ${\Bop_+}(f)$ has zeros in the
last row and column.  Furthermore ${\Aop_+}(f)$ and
${\Bop_+}(f)$ are block diagonal matrices with each block
corresponding to a $V_t\subset\pi^{-1}(V_s)$.

$\Aop_-$ and $\Bop_-$ are not homomorphism, and it is impossible to
define such a homomorphism.
\end{lemma}

\begin{proof}
For  $\rvec(n,m)\in V_s$ then,
$a_-a_+\rvec(n,m)=(C_s-\SS(s,m+1))\rvec(n,m)\ne0$ 
so $(a_-a_+)^{-1}$ is well defined via
$(a_-a_+)^{-1}\rvec(n,m)=(C_s-\SS(s,m+1))^{-1}\rvec(n,m)$. The same
is true for $(b_-b_+)^{-1}$ with
$(b_-b_+)^{-1}\rvec(n,m)=(C_s+\SS(s,m))^{-1}\rvec(n,m)$.
So $\Aop_\pm,\Bop_\pm$ are well defined.

The homomorphism of $A_+$ and $B_+$ follow from simple substitution,
as does (\ref{topch_hom}).  That $\Aop_-$ and $\Bop_-$ are not
homomorphism follows from the non existence of linear homomorphism
from the matrix algebra $M_n(\Cmpx)$ to $M_{n-1}(\Cmpx)$

\end{proof}
These maps may be considered operators for topology change. For
example in figure \ref{fig_rmline}, $J_9$ represents the
noncommutative analogue of a single connected surface as do $J_3$ and
$J_7$. We have the maps $A_+,B_+:L(V_9)\mapsto L(V_3\oplus V_7)$ this
the noncommutative analogue of a map for one surface to two. The maps 
$A_-,B_-:L(V_3\oplus V_7)\mapsto L(V_9)$ are the reverse.

%%%%%%%%%%%%%%%%%%%%%%%%%%%%%%%%%%%%%%%%%%%%%%%%%%%%%%%%%%%%%%%%%%%%%%

\subsection{Generalised Multi-topology Lattice}
\label{ch_GMTL}

Finally in this section we give an example of a further generalisation
of the multi-topology lattice $\Gspace$. This lattice no
longer refers to a specific $\rho$ but  each space $V_s$ corresponds
to a different $\rho_s$. This may have applications to the motion of
closed surfaces which, whilst they remain axially symmetric, change
shape.

There is clearly an asymmetry between the operators
${\Aop_+},{\Bop_+}$ and ${\Aop_-},{\Bop_-}$ defined by
(\ref{topch_def_ABpm}), since
${\Aop_+},{\Bop_+}$, which reflect the coalescing of
surfaces, are homomorphisms, whilst ${\Aop_-},{\Bop_-}$,
which reflect the bifurcating of a surface, are not homomorphisms. 
Thus the study of the generalised multi-topology lattice will 
enables one to study surfaces which place coalescence and
bifurcation on equal terms.

\vspace{1 em}

We define a {\bf generalised multi-topology lattice} as a Hilbert
space, $\Gspace$ which has an orthonormal basis $\{\rvec(n,m)\}$ where
the set $\{(n,m)\}$ is any subset of $\Intg\times\Intg$.  The dual of
$\rvec(n,m)$ is written $\lvec(n,m)$.  We decompose $\Gspace$ into
orthogonal subspaces
%[
\begin{align}
\Gspace &= \bigoplus_{s\in\natnum} V_s\,, 
\text{\qquad where }
V_s = \spanrm\{\rvec(n_s,m)\ , m\smin{s}\le m\le m\smax{s} \}
\end{align}
%]
The hopping operators are given by the operators
$a_\pm,b_\pm:\Gspace\mapsto\Gspace$
%[
\begin{eqnarray}
\begin{aligned}
a_+\rvec(n,m) &= 
\begin{cases}
(C_s-\SS(s,m+1))^\scrhalf\rvec(n+1,m+1) 
& \text{ if } \rvec(n,m)\in V_s
\text{ and } \rvec(n+1,m+1)\in\Gspace \\
\makebox[9 em]{0} &\text{ if } \rvec(n+1,m+1)\not\in\Gspace
\end{cases}
\\
a_-\rvec(n,m) &= 
\begin{cases}
(C_t-\SS(t,m))^\scrhalf\rvec(n-1,m-1) 
& \text{ if } \rvec(n-1,m-1)\in V_t 
\text{ for some }V_t\subset\Gspace \\
\makebox[9 em]{0} &\text{ if } \rvec(n-1,m-1)\not\in\Gspace
\end{cases}
\\
b_+\rvec(n,m) &= 
\begin{cases}
(C_s+\SS(s,m))^\scrhalf\rvec(n+1,m) 
& \text{ if }  \rvec(n,m)\in V_s
\text{ and } \rvec(n+1,m)\in\Gspace \\
\makebox[9 em]{0} &\text{ if } \rvec(n+1,m)\not\in\Gspace
\end{cases}
\\
b_-\rvec(n,m) &= 
\begin{cases}
(C_t+\SS(t,m))^\scrhalf\rvec(n-1,m) 
& \text{ if } \rvec(n-1,m)\in V_t 
\text{ for some }V_t\subset\Gspace \\
\makebox[9 em]{0} &\text{ if } \rvec(n-1,m)\not\in\Gspace
\end{cases}
\end{aligned}
\label{GMTL_abpm}
\end{eqnarray}
%]
where $\SS(s,m)\in\Real$ and $C_s\in\Real$ are constrained so that  
all the square roots are real and nonnegative. 
$\SS(s,m)$ depends on $s$ and $m$ whilst $C_s$ depends only on $s$ and
is given by $C_s=|\SS(s,m\smin{s})|=|\SS(s,m\smax{s}+1)|$

Clearly with the appropriate choice of $\SS(s,m)$ we can make $\Gspace$
correspond to one of the compromise representation of a NCSR with
over-hopping operators. A similar definition for under-hopping
operators can also be given. 

The full implications of such a lattice is being researched.

%%%%%%%%%%%%%%%%%%%%%%%%%%%%%%%%%%%%%%%%%%%%%%%%%%%%%%%%%%%%%%%%%%%%%%
%%%%%%%%%%%%%%%%%%%%%%%%%%%%%%%%%%%%%%%%%%%%%%%%%%%%%%%%%%%%%%%%%%%%%%
%%%%%%%%%%%%%%%%%%%%%%%%%%%%%%%%%%%%%%%%%%%%%%%%%%%%%%%%%%%%%%%%%%%%%%
%%%%%%%%%%%%%%%%%%%%%%%%%%%%%%%%%%%%%%%%%%%%%%%%%%%%%%%%%%%%%%%%%%%%%%

\section{Discussion}
\label{ch_disc}

This article is the result of extending the definition of tangent and
cotangent vector fields given in \cite{Gratus6} for the noncommutative
sphere, to the noncommutative surfaces of rotation given in
\cite{Gratus7}. For trivial $\rho$ we can define the algebras $\Balg$,
and $\Palg$ and use them to define such vector fields.  The results,
although more complicated, are similar to those for the case of the
noncommutative sphere.  For non-trivial $\rho$ it is not possible to
define $\Balg$ or $\Palg$. However we can still define the
multi-topology lattice and use it to examine the changes in
topology. This involved one of four possible compromises, which were
given in detail.

In order to interpret our system as a toy model for a quantised
spacetime, we need to define and interpret curvature and hence
gravity.  We indicated how one defines tangent and co-tangent vectors.
We can use the results in section \ref{ch_poi} to provide
noncommutative versions of the the metric (\ref{poi_gfh}), hodge dual
(\ref{poi_sdf}), and Laplace operator (\ref{poi_lap}).  These will all
generate problems with the choice of ordering. There is a great deal
of interest in connections and curvature \cite{Madore_curv}
\cite{Coquereaux1}. One can extend this philosophy to write the
connection $\nabla_{\wt{df}}(\wt{dh})$ and curvature
$R(\wt{df_1},\wt{df_2},\wt{df_3},df_4)$ in terms of the Poisson
structure on $\man$. However we note the r\^oles of $f$ and $h$ in
$\nabla_{\wt{df}}(\wt{dh})$ are different. Perhaps one should look at
$\nabla_X(dh)$ where $X=\ad{h}$. That is to consider the two different
types of vectors!

An alternative application of this theory is as a model of a
particle. The idea is that $\man$ is the surfaces of a particle in
space. For this we have to choose $\rho$, an element $\Hamil\in\Ralg$
such that $\Hamil=\Hamil^\dagger$ and a value of $\varepsilon$. Once
these are fixed we can calculate the spectrum of $\Hamil$. If we
interpret $\Hamil$ as a Hamiltonian, then we could call the spectrum
of $\Hamil$ a mass spectrum. The main problem with this is choosing
$\rho$ and $\Hamil$ since there is not much physics we can use to
guide us. If $\rho$ is not trivial and we choose $\Hamil$ to contain
the hopping operators then this could be used as a model of
interacting particles.

%%%%%%%%%%%%%%%%%%%%%%%%%%%%%%%%%%%%%%%%%%%%%%%%%%%%%%%%%%%%%%%%%%%%%%

\subsection*{Acknowledgement}

The author would like to thank Robin Tucker and Marianne
Karlsen for their suggestions and help in the preparation of this
article, and the physics department of Lancaster University for there
facilities.

%%%%%%%%%%%%%%%%%%%%%%%%%%%%%%%%%%%%%%%%%%%%%%%%%%%%%%%%%%%%%%%%%%%%%%
%%%%%%%%%%%%%%%%%%%%%%%%%%%%%%%%%%%%%%%%%%%%%%%%%%%%%%%%%%%%%%%%%%%%%%
%%%%%%%%%%%%%%%%%%%%%%%%%%%%%%%%%%%%%%%%%%%%%%%%%%%%%%%%%%%%%%%%%%%%%%
%%%%%%%%%%%%%%%%%%%%%%%%%%%%%%%%%%%%%%%%%%%%%%%%%%%%%%%%%%%%%%%%%%%%%%

%%%%%%%%%%%%%%%%%%%%%%%%%%%%%%%%%%%%%%%%%%%%%%%%%%%%%%%%%%%%%%%%%%%%%%
%%%%%%%%%%%%%%%%%%%%%%%%%%%%%%%%%%%%%%%%%%%%%%%%%%%%%%%%%%%%%%%%%%%%%%
%%%%%%%%%%%%%%%%%%%%%%%%%%%%%%%%%%%%%%%%%%%%%%%%%%%%%%%%%%%%%%%%%%%%%%

\newpage

\begin{table}

\abovedisplayshortskip 0pt
\belowdisplayshortskip 0pt
\abovedisplayskip 0pt
\belowdisplayskip 0pt
\topsep 0 pt

\hspace{-1cm}
\begin{tabular}{|p{2 truecm}|p{7.5 truecm}|p{7.5 truecm}|}
\hline
Name &
$\Calg=\Calg(\rho,\varepsilon)$ &
$\Ralg=\Ralg(\rho,\varepsilon,R)$ 
\\
\hline
Parameters &
$\rho\in C^1(\Real\mapsto\Real)$,
$\varepsilon\in\Real$ 
&
$\rho\in C^1(\Real\mapsto\Real)$,
$R,\varepsilon\in\Real$ 
\\
Generators &
$X_+,X_-$, and $f(X_0)$ where $f\in\Xdom$ 
&
$X_+,X_-$, and $f(X_0)$ where $f\in\Xdom$ 
\\
\hline
\begin{tabular}[c]{@{}l@{}}
General \\
\ normal \\
\ ordered \\
\ element 
\end{tabular}
&
\begin{minipage}[c]{7.5 truecm}
%[
\begin{eqnarray}
f = 
\sum_{r=0}^{r_{\max}}
\sum_{s=0}^{s_{\max}}
X_+^r X_-^s f_{rs}(X_0)
\label{rev_Ca_ge}
\\
\text{ with } 
f_{rs}\in\Xdom
\qquad\qquad
\notag
\end{eqnarray}
%]
\end{minipage}
&
\begin{minipage}[c]{7.5 truecm}
%[
\addtocounter{equation}{2}
\begin{eqnarray}
&&
f = 
\sum_{r=0}^{r_{\max}}
X_\pm^r f_{r}(X_0) 
\label{rev_Ra_ge}
\\
&&
\text{ with }
f_{r}\in\Xdom
\qquad\qquad
\notag
\end{eqnarray}
%]
\end{minipage}
\\
\hline
\begin{tabular}[c]{@{}l@{}}
Quotient \\
\ equations
\end{tabular}
&
\begin{minipage}[c]{7.5 truecm}
%[
\addtocounter{equation}{-3}
\begin{eqnarray}
f(X_0)X_\pm &=& X_\pm f(X_0\pm\varepsilon)
\label{rev_Ca_f}
\\{}
[X_+,X_-] &=& \rho(X_0+\varepsilon)-\rho(X_0)
\label{rev_Ca_XX}
\end{eqnarray}
%]
\end{minipage}
&
\begin{minipage}[c]{7.5 truecm}
%[
\addtocounter{equation}{1}
\begin{eqnarray}
f(X_0)X_\pm &=& X_\pm f(X_0\pm\varepsilon)
\label{rev_Ra_f}
\\
X_+X_- &=& \rho(R)-\rho(X_0)
\label{rev_Ra_XX}
\end{eqnarray}
%]
\end{minipage}
\\
\hline
%%%%%%%%%%%%%%%%%%%%%%%%%%%%%%%%%%%%%%%%%%%%%%%%%%%%%%%%%%
\multicolumn{3}{l}{\strut}
\\
\hline
Name &
$\Nalg=\Nalg(\rho,\varepsilon)$ &
$\NCalg=\NCalg(\rho,\varepsilon,R)$ 
\\
\hline
Parameters &
$\rho\in C^1(\Real\mapsto\Real)$,
$\varepsilon\in\Real$ 
&
$\rho\in C^1(\Real\mapsto\Real)$,
$\varepsilon\in\Real$ 
\\
Generators &
$X_+,X_-$, and $f(X_0,N_0)$ where $f\in\NXdom$ 
&
$X_+,X_-$, and $f(X_0,N_0)$ where $f\in\NXdom$ 
\\
\hline
\begin{tabular}[c]{@{}l@{}}
General \\
\ normal \\
\ ordered \\
\ element 
\end{tabular}
&
\begin{minipage}[c]{7.5 truecm}
%[
\begin{eqnarray}
&&
f = 
\sum_{r=-r_{\max}}^{r_{\max}} 
X_\pm^r f_{r}(X_0,N_0) 
\label{rev_Na_ge}
\\
&&
\text{ with }
f_{r}\in\NXdom
\notag
\end{eqnarray}
%]
\end{minipage}
&
\begin{minipage}[c]{7.5 truecm}
%[
\addtocounter{equation}{3}
\begin{eqnarray}
&&
f = 
\sum_{r=0}^{r_{\max}}
\sum_{s=0}^{s_{\max}}
X_+^r X_-^s f_{rs}(X_0,N_0) 
\label{rev_NCa_ge}
\\
&&
\text{ with }
f_{rs}\in\NXdom
\qquad\qquad
\notag
\end{eqnarray}
%]
\end{minipage}
\\
\hline
\begin{tabular}[c]{@{}l@{}}
Quotient \\
\ equations
\end{tabular}
&
\begin{minipage}[c]{7.5 truecm}
%[
\addtocounter{equation}{-4}
\begin{eqnarray}
f(X_0,N_0)X_\pm &=& X_\pm f(X_0\pm\varepsilon,N_0)
\qquad
\label{rev_Na_f}
\\
X_+X_- &=& \rho(N_0)-\rho(X_0)
\label{rev_Na_XX}
\\{}
[X_0,N_0] &=& 0
\label{rev_Na_XN}
\end{eqnarray}
%]
\end{minipage}
&
\begin{minipage}[c]{7.5 truecm}
%[
\addtocounter{equation}{1}
\begin{eqnarray}
f(X_0,N_0)X_\pm &=& X_\pm f(X_0\pm\varepsilon,N_0)
\qquad
\label{rev_NCa_f}
\\{}
[X_+,X_-] &=& \rho(X_0+\varepsilon)-\rho(X_0)
\label{rev_NCa_XX}
\\{}
[X_0,N_0] &=& 0
\label{rev_NCa_XN}
\end{eqnarray}
%]
\end{minipage}
\\
\hline
\end{tabular}
\caption{Properties of the algebras $\Calg,\Ralg,\Nalg,\NCalg$}
\label{rev_tbl}
\end{table}
%%%%%%%%%%%%%%%%%%%%%%%%%%%%%%%%%%%%%%

\newpage

\begin{figure}
\setlength{\unitlength}{60 mm}
\begin{picture}(2.5,2.5)
\put(0,1){\epsfxsize\unitlength
\epsffile{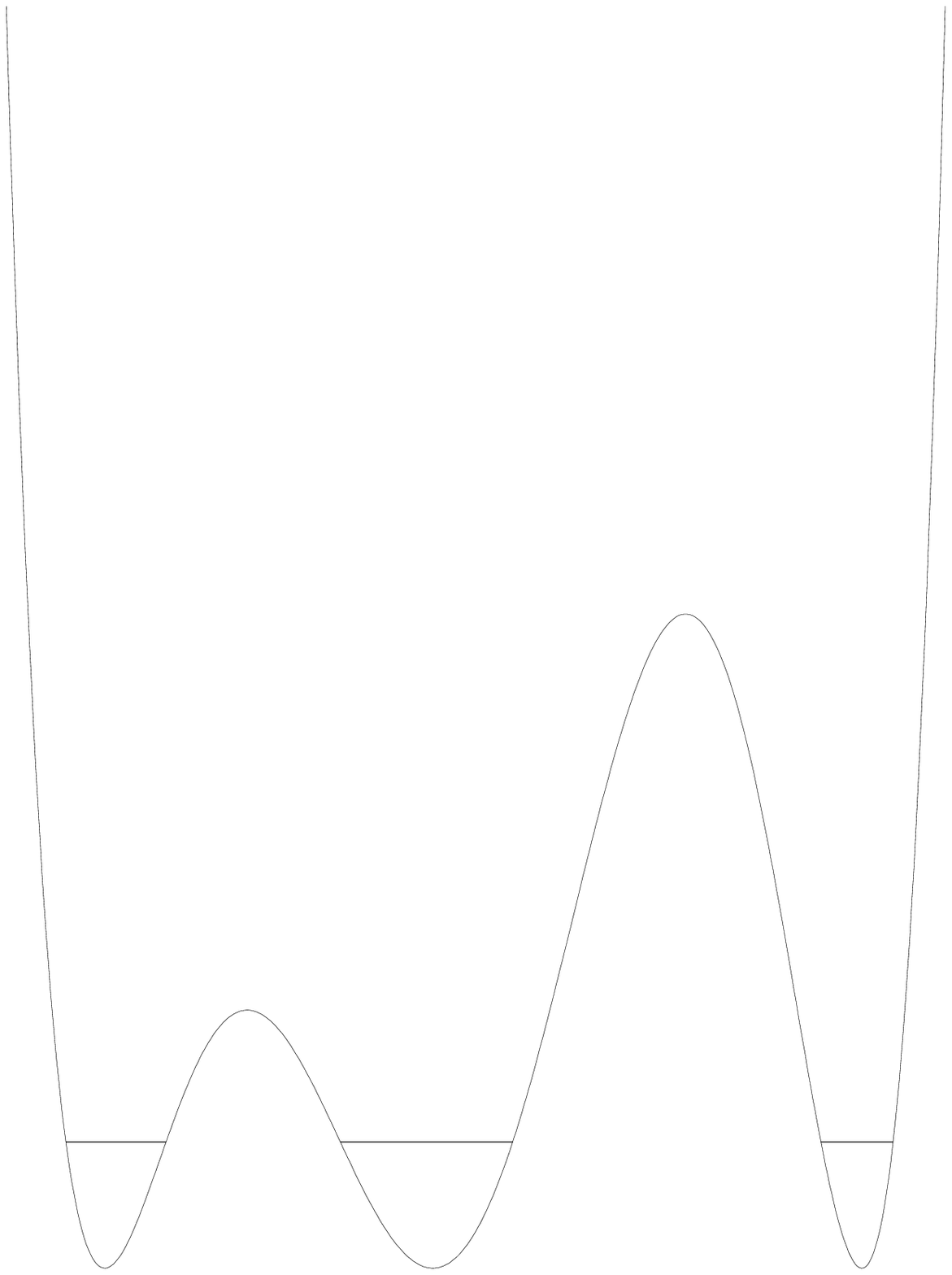}}
\put(1,1){\makebox(0,0)[bl]{\rotatebox{270}{\epsfxsize\unitlength
\epsffile{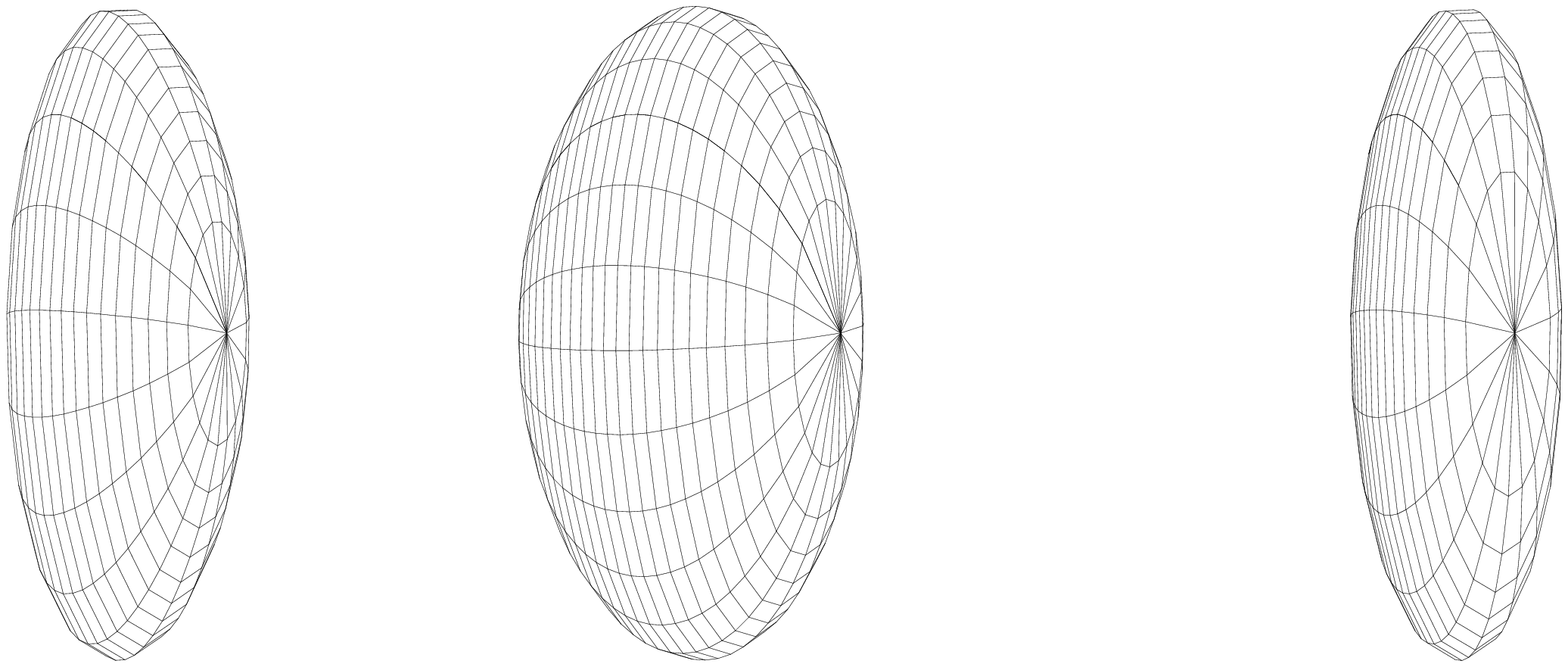}}}}
\put(0,0){\epsfxsize\unitlength
\epsffile{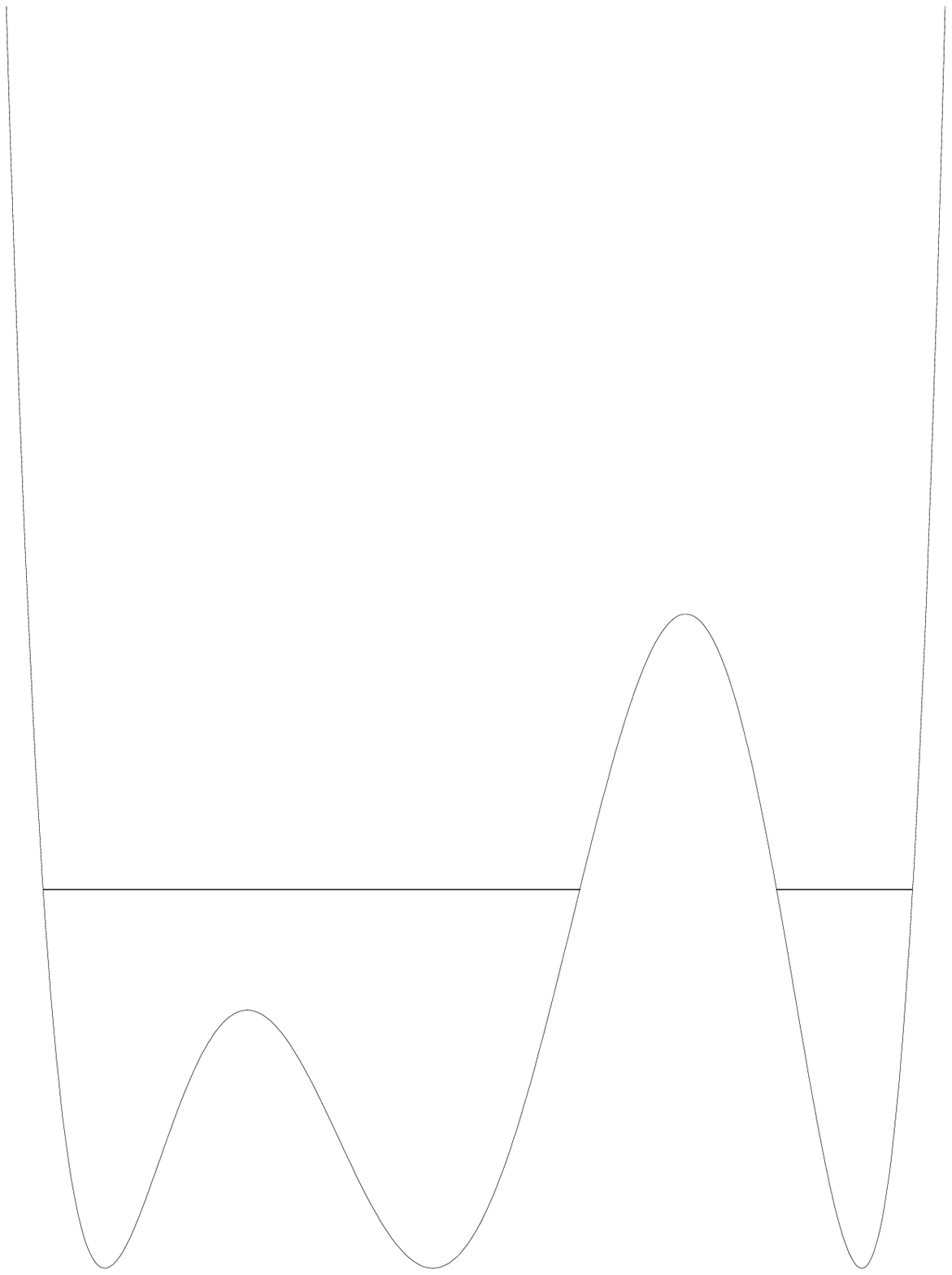}}
\put(1,0){\makebox(0,0)[bl]{\rotatebox{270}{\epsfxsize\unitlength
\epsffile{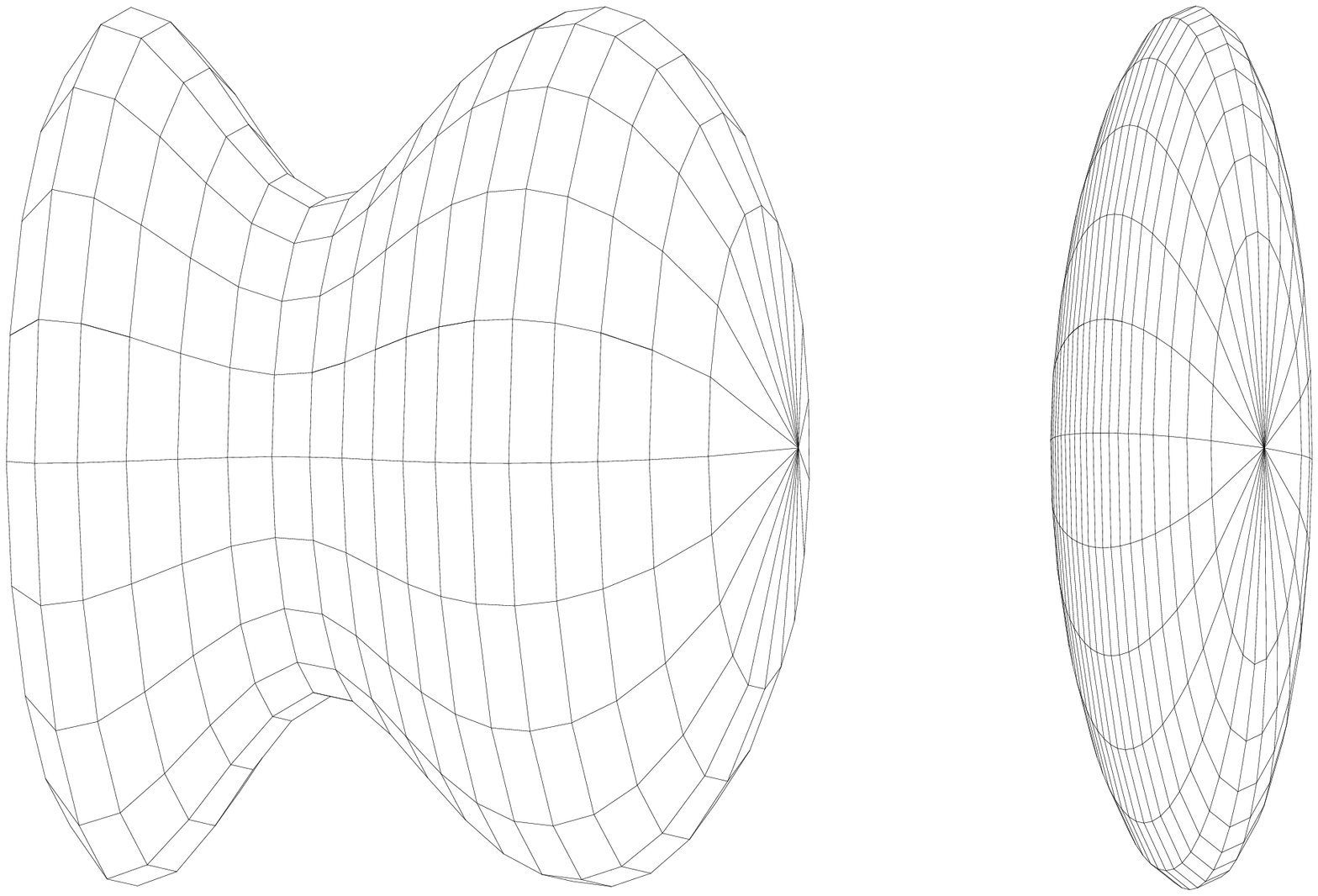}}}}
\put(0,0){\setlength{\unitlength}{6 mm}
\begin{picture}(25,25)
%%\graphpaper[1](0,0)(25,25)
\put(1.5,13.5){\makebox(0,0)[b]{$I_1$}}
\put(4.0,13.5){\makebox(0,0)[b]{$I_2$}}
\put(7.5,13.5){\makebox(0,0)[b]{$I_3$}}
\put(12,18){\makebox(0,0)[b]{$\man_1$}}
\put(16.5,18){\makebox(0,0)[b]{$\man_2$}}
\put(23,18){\makebox(0,0)[b]{$\man_3$}}
\put(3,5.6){\makebox(0,0)[b]{$I_4$}}
\put(7.5,5.6){\makebox(0,0)[b]{$I_5$}}
\put(15,9){\makebox(0,0)[b]{$\man_4$}}
\put(22,9.5){\makebox(0,0)[b]{$\man_5$}}
\end{picture}}
\end{picture}
\caption{Surface of rotation corresponding to the same $\rho$ for
different values of $R$}
\label{fig_top}
\vspace{10em}
\end{figure}

\newpage

\begin{figure}
\begin{center}
\setlength{\unitlength}{40 mm}
\begin{picture}(2,1.2)(0,-0.2)
\put(0,0){\epsfysize\unitlength \epsffile{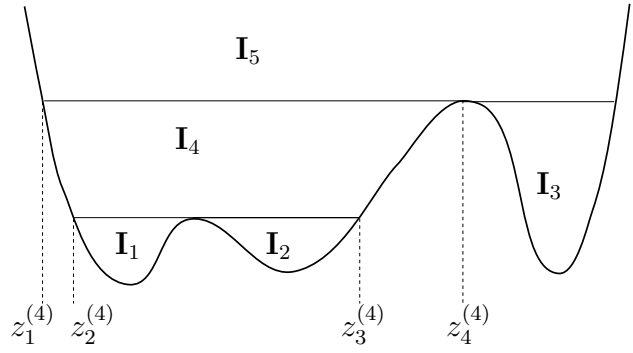}}
\put(0.3,0.15){\makebox(0,0)[bl]{$\topI{1}$}}
\put(0.8,0.15){\makebox(0,0)[bl]{$\topI{2}$}}
\put(1.7,0.35){\makebox(0,0)[bl]{$\topI{3}$}}
\put(0.5,0.5){\makebox(0,0)[bl]{$\topI{4}$}}
\put(0.7,0.8){\makebox(0,0)[bl]{$\topI{5}$}}
\put(-0.05,0){\makebox(0,0)[tl]{$z^{(4)}_1$}}
\put(0.15,0){\makebox(0,0)[tl]{$z^{(4)}_2$}}
\put(1.05,0){\makebox(0,0)[tl]{$z^{(4)}_3$}}
\put(1.4,0){\makebox(0,0)[tl]{$z^{(4)}_4$}}
\end{picture}
\end{center}
\caption{example of topology intervals}
\label{fig_topi}
\vspace{10em}
\end{figure}

\newpage

\newpage

\begin{figure}[t]
\begin{center}
\epsfysize 5 cm 
\centerline{\epsffile{nontriv.eps}} 
\end{center}
\setlength{\unitlength}{10 mm}
\begin{picture}(17.5,5.7)(-1,0.8)
\multiput(0,6)(2,0){9}{$\rvec(6,{\arabic{nn}})$\addtocounter{nn}{1}}
\put(0,0){\setcounter{nn}{0}}
\multiput(1,5)(2,0){8}{$\rvec(5,{\arabic{nn}})$\addtocounter{nn}{1}}
\put(0,0){\setcounter{nn}{0}}
\multiput(2,4)(2,0){7}{$\rvec(4,{\arabic{nn}})$\addtocounter{nn}{1}}
\put(0,0){\setcounter{nn}{0}}
\multiput(3,3)(2,0){2}{$\rvec(3,{\arabic{nn}})$\addtocounter{nn}{1}}
\put(4,2){$\rvec(2,0)$}
\put(0,0){\setcounter{nn}{3}}
\multiput(9,3)(2,0){3}{$\rvec(3,{\arabic{nn}})$\addtocounter{nn}{1}}
\put(0,0){\setcounter{nn}{3}}
\multiput(10,2)(2,0){2}{$\rvec(2,{\arabic{nn}})$\addtocounter{nn}{1}}
\put(11,1){$\rvec(1,3)$}
\put(9.3,5.3){\vector(1,1){0.6}}
\put(9.3,5){\vector(1,-1){0.6}}
\put(9,5.3){\vector(-1,1){0.6}}
\put(9,5){\vector(-1,-1){0.6}}
\put(8.3,5.4){$b_+$}
\put(8.3,4.8){$a_-$}
\put(9.7,5.4){$a_+$}
\put(9.7,4.8){$b_-$}
\put(-1,6){$V_8\,\to$}
\put(0,5){$V_7\,\to$}
\put(1,4){$V_6\,\to$}
\put(2,3){$V_2\,\to$}
\put(3,2){$V_1\,\to$}
\put(8,3){$V_5\,\to$}
\put(9,2){$V_4\,\to$}
\put(10,1){$V_3\,\to$}
\put(11.6,5.1){\vector( 1,0){1.2}}
\put(12.8,5.0){\vector(-1,0){1.2}}
\put(12,5.2){$X_+$}
\put(12,4.7){$X_-$}
\vspace{10em}
\end{picture}

\begin{center}[t]
\begin{tabular}{|c|c|l|}
\hline
space & dim & basis \\
\hline
$V_1$ & 1 & $\{\rvec(2,0)\}$ \\
$V_2$ & 2 & $\{\rvec(3,0),\rvec(3,1)\}$ \\
$V_3$ & 1 & $\{\rvec(1,3)\}$ \\
$V_4$ & 2 & $\{\rvec(2,3),\rvec(2,4)\}$ \\
\hline
\end{tabular}
\hspace{4 em}
\begin{tabular}{|c|c|l|}
\hline
space & dim & basis \\
\hline
$V_5$ & 3 & $\{\rvec(3,3),\rvec(3,4),\rvec(3,5)\}$ \\
$V_6$ & 7 & $\{\rvec(4,0),\ldots,\rvec(4,6)\}$ \\
$V_7$ & 8 & $\{\rvec(5,0),\ldots,\rvec(5,7)\}$ \\
$V_8$ & 9 & $\{\rvec(6,0),\ldots,\rvec(6,8)\}$ \\
\hline
\end{tabular}
\end{center}
\caption{The First Few Representations and the multi-topology lattice
for a non-trivial $\rho$}
\label{fig_notriv}
\vspace{10em}
\end{figure}

\newpage

\newpage

\newpage

\newpage

\newpage

\newpage

\end{document}